\documentclass[12pt]{amsart}
\usepackage{amsmath,amssymb,amsthm}
\usepackage{setspace}
\usepackage{enumerate}
\usepackage{tikz}

\usepackage{amsmath,amssymb,amsthm}
\usepackage{comment}
\newcommand{\R}{\mathbb{R}}
\newcommand{\Nat}{\mathbb{N}}
\newcommand{\Q}{\mathbb{Q}}
\newcommand{\B}{\mathcal{B}}

\newcommand{\restrict}{\upharpoonright}

\newcommand{\C}{2^\omega}

\newcommand{\dom}{\operatorname{dom}}
\newcommand{\kO}{\mathcal O}
\newcommand{\eps}{\varepsilon}
\newcommand{\concat}{^\smallfrown}
\newcommand{\Ord}{\operatorname{Ord}}

\newcommand{\leqm}{\leq_{\mathbf m}}
\newcommand{\geqm}{\geq_{\mathbf m}}
\newcommand{\leqtt}{\leq_{\mathbf{tt}}}

\newcommand{\leqtto}{\leq_{\mathbf{tt1}}}

\newcommand{\leqT}{\leq_{\mathbf T}}
\newcommand{\geqT}{\geq_{\mathbf T}}
\newcommand{\Baire}{\omega^\omega}

\newtheorem{thm}{Theorem}
\newtheorem{defn}[thm]{Definition}
\newtheorem{prop}[thm]{Proposition}
\newtheorem{cor}[thm]{Corollary}
\newtheorem{lemma}[thm]{Lemma}
\newtheorem{question}[thm]{Question}

\newcommand{\BS}[1]{\mathbf{{\Sigma}}^{0}_{#1}}

\newcommand{\BP}[1]{\mathbf{{\Pi}}^{0}_{#1}}

\title[Three topological reducibilities]{Three topological reducibilities for discontinuous functions}

\author[A.\ R.\ Day]{Adam R. Day}
\address{School of Mathematics and Statistics\\
  Victoria University of Wellington\\
  Wellington, New Zealand}
\email{adam.day@vuw.ac.nz}

\author[R.\ Downey]{Rod Downey}
\address{School of Mathematics and Statistics\\
  Victoria University of Wellington\\
  Wellington, New Zealand}
\email{rod.downey@vuw.ac.nz}

\author[L.\ B.\ Westrick]{Linda Brown Westrick}
\address{Department of Mathematics\\
Penn State University\\
University Park, Pennsylvania U.S.A.}
\email{westrick@psu.edu}

\thanks{This work was supported in part by the Marsden Fund of New Zealand.  Westrick was supported in part by Noam Greenberg's Rutherford Discovery Fellowship as a postdoctoral fellow.}

\begin{document}

\begin{abstract}
We define a family of three related reducibilities, $\leqT$, $\leqtt$ 
and $\leqm$, for 
arbitrary functions
$f,g:X\rightarrow\mathbb R$, where $X$ 
is a compact separable metric space.
The $\equiv_{\mathbf T}$-equivalence classes 
mostly coincide with the proper Baire classes.
We show that certain $\alpha$-jump functions 
$j_\alpha:2^\omega\rightarrow \R$ 
are $\leqm$-minimal in their Baire class.
Within the Baire 1 functions, we completely
characterize the degree structure associated 
to $\leqtt$ and $\leqm$, finding an exact match to
 the $\alpha$ hierarchy introduced by
Bourgain \cite{Bourgain}
and analyzed in Kechris-Louveau \cite{KechrisLouveau1990}.

\end{abstract}

\maketitle

\section{Introduction}
\subsection{Reducibilities}
Computability theory seeks to understand the effective content of 
mathematics. 
Ever since its beginnings in the work of G\"odel, Turing, Post, Kleene,
Church and others, the idea of a reduction has been a central notion
in this area. Turing \cite{Turing1939} formalized  we now call \emph{Turing reducibility}
which can be viewed as the most general way of allowing computation
of one set of natural numbers from another using oracle queries.

In the last 60 years, we have seen the introduction of a large number of 
reducibilities $A\le B$, reflecting different access mechanisms 
for the computation of $A$ from $B$. 
Different oracle access mechanisms give different 
equivalence classes calibrating computation.
The measure of the efficacy of such reductions
is the extent to which 
\begin{itemize}
\item[(i)] they give insight into computation, and 
\item[(ii)] they are useful in mathematics.
\end{itemize}
Examples of (ii) above,
include the use of 
polynomial time reductions to enable the  theory of $NP$-completeness, but also
include the use of $\Pi_1^1$-completeness to demonstrate that 
classical isomorphism problems like the classification of countable 
abelian groups cannot have reasonable invariants (Downey-Montalb{\'a}n \cite{DowneyMontalban2008}),
Ziegler reducibility to classify algebraic closures of finitely presented 
groups (see e.g. Higman-Scott \cite{HigmanScott1988}), 
truth-table reducibility to analyze algorithmic randomness for continuous 
measures (Reimann-Slaman \cite{ReimannSlamanTA}), and enumeration reducibility 
for the relativised Higman embedding theorem (see \cite{HigmanScott1988}). 
There are many other examples.

\subsection{Reducibilities in type II computation}
The narrative above really only refers to notions of relative computability
for infinite bit sequences (or objects, such as real numbers, 
which can be coded by such sequences).  That is, the 
objects whose information content is being compared have 
function type $A:\omega\rightarrow \omega$ or similar.

What if instead we wanted to compare the information content 
of functions $f:[0,1]\rightarrow\mathbb R$?
The collection $\mathcal F([0,1])$ of all such functions has 
cardinality greater than the continuum, so it is not possible 
to use infinite bit sequences to code all these objects.  In the 
next section we will say a bit more about 
some approaches to the problem
of relative computability for higher type objects, the most 
prominent of which is the Weihrauch computable reducibility 
framework.

In this paper, we introduce and analyse three  notions of reduction for 
$\mathcal F(X)$, where $X$ is a compact Polish space. 
Two of our notions are completely new and one has had little
previous attention. We argue that that they meet the criteria 
(i) and (ii) above, and provide computational insight into the 
hierarchies previously introduced in classical analysis for the classification
of the Baire classes of functions\footnote{We define these terms in Section
\ref{preliminaries}.}.

We first concentrate upon what we define  to be 
$f\le_{\mathbf T} g$. This reduction is interpreted to mean 
that $f$ is continuously Weihrauch reducible 
to the \emph{parallelization} of $g$. In the next section, 
we define what we mean by this, and  argue that this is the 
most natural (continuous) analog of 
Turing reducibility for higher type objects. 
We introduce the new notions of
$f\leq_{\mathbf{tt}} g$ and $f\leq_{\mathbf m} g$
by restricting the oracle use of the functionals
in the Weihrauch reduction in an appropriate way described
in Section \ref{sec:Definition_of_tt_and_m_reducibility}.

It seems to be folklore that  the $\leq_{\mathbf T}$
degrees of the Baire functions are linearly ordered, and 
these degrees correspond to the proper Baire classes. 
Our main 
results concern the $\leq_{\mathbf m}$ and $\leq_{\mathbf{tt}}$ 
degrees.  We show that the $\alpha$th
jump operator\footnote{We will define $j_\alpha$ later, but for example 
$j_1:2^\omega \rightarrow \mathbb R$ is $j_1(X) := \sum_{i \in X'} 2^{-(i+1)}$.}
 $j_\alpha$ is $\leq_{\mathbf m}$-minimal in its Baire class. 
\begin{thm}
If a Baire function $f$ is not Baire $\alpha$, then $f \geqm j_{\alpha +1}$.
\end{thm}

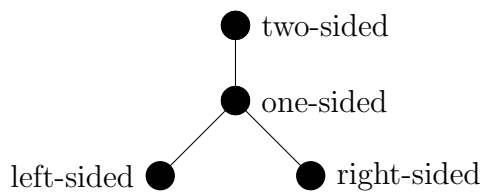
\begin{figure}
\begin{center}
\begin{tikzpicture}
\draw node[fill,circle, label={right:one-sided}] at (0,3) {};
\draw node[fill,circle, label={right:two-sided}] at (0,4) {};
\draw node[fill,circle, label={left:left-sided}] at (-1,2) {};
\draw node[fill,circle, label ={right:right-sided}] at (1,2) {};
\draw node[fill,circle] at (0,3) {};
\draw node[fill,circle] at (0,4) {};
\draw (1,2) -- (0,3) -- (0,4);
\draw (-1,2) -- (0,3);
\end{tikzpicture}
\end{center}
\caption{The $\leqm$ degrees of functions $f$ with $|f|_\alpha$ a successor.}\label{fig:m-degrees}
\end{figure}

Then we restrict attention to the Baire 1 functions.  
In \cite{KechrisLouveau1990}, Kechris and Louveau  consider 
three ranking functions 
$\alpha$, $\beta$ and $\gamma$,  
which take Baire 1 functions to countable ordinals. 
These ranks are especially robust at levels of the form
$\omega^\xi$. 
Letting $\xi(f)$ denote the least $\xi$ such that 
$\alpha(f) \leq \omega^\xi$, in our main theorem 
we characterize the $\leqm$ and $\leqtt$ 
degrees of the Baire 1 functions as follows.
\begin{thm}\label{thm:T2} For $f$ and $g$ discontinuous Baire 1 functions,
\begin{enumerate}
\item $f\leqtt g$ if and only if $\xi(f) \leq \xi(g)$.
\item If $|f|_\alpha < |g|_\alpha$, then $f\leqm g$.
\item If $\nu$ is a limit ordinal, $\{f : |f|_\alpha = \nu\}$ is 
an $\leqm$-degree.
\item If $\nu$ is a successor, $\{f: |f|_\alpha = \nu\}$ 
contains exactly four $\leqm$-degrees arranged as in Figure \ref{fig:m-degrees}.
\end{enumerate}
\end{thm}

\begin{figure}
\begin{center}
\begin{tikzpicture}
\draw node[fill,circle, label={right:constant}] at (0,0) {};
\draw node[fill,circle, label={right:continuous}] at (0,1) {};
\draw node[fill,circle, label={left:lower-semi-continuous}] at (-1,2) {};
\draw node[fill,circle, label ={right:upper-semi-continuous}] at (1,2) {};
\draw node[fill,circle] at (0,3) {};
\draw node[fill,circle] at (0,4) {};
\draw (-5.5,1.7) rectangle (5.5,4.5);
\draw (-5.5,1.3) rectangle (5.5, -.3);
\draw (0,0) -- (0,1) -- (1,2) -- (0,3) -- (0,4);
\draw (0,1) -- (-1,2) -- (0,3);
\node at (-3.5,.5) {$|f|_\alpha = 1$};
\node at (-3.5,3.5) {$|f|_\alpha = 2$};
\end{tikzpicture}
\end{center}
\caption{The smallest $\leqm$ degrees are recognizable classes.}\label{fig:recognizable}
\end{figure}
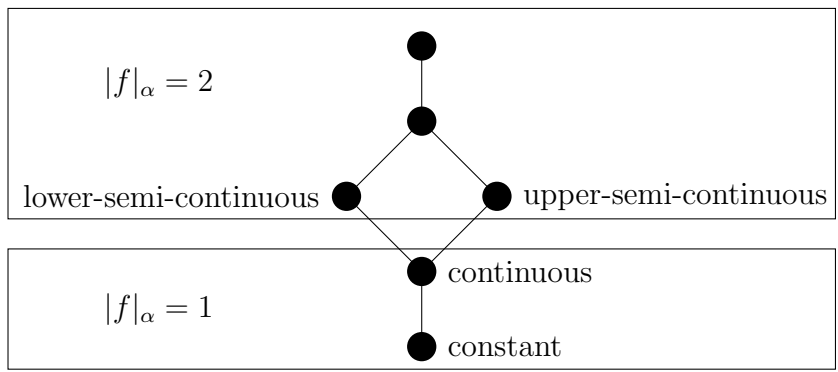

The smallest $\leqm$-degrees are recognizable classes: constant functions,
continuous functions, upper semi-continuous functions, and 
lower semi-continuous functions.  See Figure \ref{fig:recognizable}.

The authors would like to thank Vasco Brattka, Takayuki Kihara, 
Antonio Montalb\'an, Arno Pauly 
and Dan Turetsky for many useful discussions on these and related topics.

\section{Motivations in defining $\le_{\mathbf T}$}
\subsection{Weihrauch/computable  reducibility}
Suppose we want to define $f\le g$ for functions $f$ and $g$, with the 
meaning that $g$ can compute $f$. The search for a natural notion of 
$f\le g$ leads directly to Weihrauch reducibility.  For $A,B \in 2^\omega$,
it is clear what it means to ``know'' $A$.  An algorithm or oracle 
knows $A$ if, given input $n$, it outputs $A(n)$.  Accordingly, 
a computation of $A$ from $B$ is a algorithm which can 
answer these questions about $A$ when given query access to an 
oracle for $B$.  So, what kinds of questions should we be able 
to answer if we claim to ``know'' $f:[0,1]\rightarrow \R$?  At a minimum,
an oracle for $f$ ought to be able to 
produce $f(x)$ when given input $x$.  We take this ability 
as the defining feature of an oracle for $f$.

Now, what should it mean
for an algorithm to have query access to an oracle for $g$?
Clearly, given input $x$, the algorithm should be able to
pass it through and 
query $g(x)$.  If $g(x)$ were the only permitted query, 
the algorithm could not really be said to have access to an 
oracle for all of $g$, so we should allow some other queries 
as well.  For example, one would hope for a theory in which
 the functions $x\mapsto f(x)$ and $x\mapsto f(x+c)$ 
 always compute each other, where $c$ is a computable 
 real.  Generalizing this idea, an algorithm with 
 query access to $g$ should be able to ask about $g(y)$ 
 for any $y\leq_T x$.  Therefore, the notion of Weihrauch 
 reducibility is a natural starting candidate 
 for a notion of $f \le g$.  For simplicity 
 of this narrative we horribly abuse some notation 
 in the definition given below.
 We refer the reader to Section \ref{reps} for the 
 full definition.

\begin{defn} Let $f, g:[0,1]\rightarrow \R$.
\begin{itemize}
\item[(a)]
 $f$ is called \emph{computably/Weihrauch reducible to 
$g$}, $f\le_W g$, if there are partial computable functions 
$H,K:2^\omega\rightarrow 2^\omega$ with 
$$f(x)=H(x,g(K(x))).$$
\item[(b)] If we replace ``computable''
by continuous in the definition above, we refer to this as 
\emph{continuous} Weihrauch reducibility, and write 
$f\le_W^c g$.
\end{itemize}
\end{defn}

The abuse of notation lies in the fact that $H$ and $K$ 
manipulate names or codes for reals rather than the reals
themselves.\footnote{A subtle but important point is 
that reals have multiple names, and $K$ is 
not required to have consistent behavior on  
two different names for the same real.}

The name ``Weihrauch reducibility'' was coined by 
Brattka and Gherardi \cite{BrattkaGherardi2011}, whereas earlier Weihrauch had called 
it computable reducibility. 

Computable/continuous Weihrauch
reducibility has been
studied in \cite{brattka2005}, \cite{mylatz2006}
and \cite{pauly2010}.
Brattka \cite{brattka2005}
proved effective versions of classical
theorems linking the Borel and Baire hierarchies
using this reducibility.

\subsection{Parallelized Weihrauch reducibilty}

The above account seems to miss the feature of 
 ordinary $\leq_T$ computation in which the algorithm 
may use the oracle repeatedly and interactively.  
We would not like to 
limit the reduction algorithm to a single use of the 
$g$ oracle.

However, if the algorithm had access to all of $g(x)$ 
based on its first query, it would be able to feed this 
back into the $g$ oracle, obtaining $g(g(x))$ 
and in general the sequence of $g^{(n)}(x)$.
And if we accept 
some algorithm is 
uniformly producing the sequence $g^{(n)}(x)$, 
it could be simultaneously engaged in writing down 
a summarizing output $g^{(\omega)}(x)$, where 
$g^{(\omega)}(x)$ is for example defined as
$\bigoplus_n g^{(n)}(x).$\footnote{Imagine 
for the purposes of this hypothetical that $g$ is an operator 
on $2^\omega$, so that a joining operation $\oplus$ is available 
to us; a similar situation could be concocted 
for operators on the unit interval.}  So we are led to accept 
$g^{(\omega)} \le g$.
If we accept this, and also wish our notion to be transitive, 
we must accept $g^{(\omega+1)} \le g$, 
otherwise transitivity will be violated in the sequence 
$g^{(\omega+1)} \le (g^{(\omega)} \oplus g) \le g$.  
In the end, we are forced to say $g$ computes all 
its iterates up to $\omega_1^{ck}$.  The notion just described, 
complete with all the transfinite iteration, was studied by 
Kleene \cite{Kleene1959}.  However, this reducibility is coarser than we want
(for example, we would not want the 
jump operator on $2^\omega$ to be able to compute 
the double-jump operator) 
and so we choose to go by another route.

Suppose instead we make the following seemingly minor 
adjustment to our concept of what an oracle for $g$ should do. 
Instead of querying with an input $x$, we query with a
pair $(x,\eps)$, where $\eps \in \Q^+$.  Instead of returning 
the entire $g(x)$, 
the oracle returns some $p \in \Q$ with $|g(x) -p|<\eps$.
Now an algorithm which on input $x$ has made finitely 
many queries to $g$ has only acquired a finite amount of 
new information, so its future queries are still restricted 
to those $y$ with $y \leq_T x$.  This breaks the cycle above.
In order to get more and more precision on $f(x)$,
such an algorithm may query $g(y)$ for many different 
values of $y$.  But there are at most countably many 
queries to $g$ associated to the computation of a single $f(x)$.
Therefore, we can naturally express the kind of 
reducibility described above in the 
Weihrauch framework: $f \le g$ could mean
$f \leq_{W} \hat g$, where 
$\hat g:X^\omega \rightarrow Y^\omega$ is the \emph{parallelization} of $g$,
defined by applying $g$ componentwise.

\subsection{What is a single bit of information about $f$?}

Accepting parallelized Weihrauch reducibility as a higher-type
notion of $\leq_T$, what should $\leq_{tt}$ and $\leq_m$ be?
It is particularly informative to consider $\leq_m$.  In 
classical computability theory, $A\leq_m B$ means that 
there is an algorithm which, on input $n$, outputs $m$ 
such that $A(n) = B(m)$.  For us the important features are:
\begin{enumerate}
\item The oracle's response is accepted unchanged as the output, and 
\item The question is a yes/no question.
\end{enumerate}
Allowing a more demanding question (such as ``approximate $f(x)$ 
to within $\eps$'') seems unfair, ruling out $m$-computations 
between functions of disjoint ranges that are otherwise
computationally identical.  (Notions of $m$-reducibility 
without point (2) have been considered however, 
for example by Hertling \cite{Hertling1993}, Pauly \cite{pauly2010} and Carroy \cite{carroy2013}.)

Our previous decision on how to finitize the oracle 
was, upon reflection, rather arbitrary.  We could 
restrict ourselves to yes/no questions with the following 
convention about oracles, and still end up with a $\leq_T$ notion 
equivalent to
parallelized Weihrauch reducibility.
An oracle for $f$ accepts as input a triple $(x,p,\eps)$,
with $p \in \Q$ and $\eps \in \Q^+$, and $\eps$-approximately 
answers the question ``is $f(x) < p$''?  The exact 
version of this question would be too precise for a computable 
procedure, so we accept any answer as correct 
if $|f(x) - p| < \eps$.  Now that each query to the oracle yields 
exactly one bit of information, we can define $\leq_m$ 
and $\leq_{tt}$ for the higher type objects
by placing corresponding restrictions 
on the oracle use.  We do this in Section \ref{sec:Definition_of_tt_and_m_reducibility}.

\subsection{Parameters} Another natural question we might 
ask ourselves is ``what parameters would be 
reasonable for such reductions''?  For reductions between
objects of type $A:\omega\rightarrow \omega$, we usually 
allow integer parameters in computation procedures.
Therefore, for reductions between objects of type 
$f:[0,1]\rightarrow \R$, perhaps we should allow real 
parameters.  We take this approach, which has a 
substantial simplifying effect.  Every continuous function 
is computable relative to a real parameter, so Weihrauch 
computability relative to a real parameter is the same as 
continuous Weihrauch reducibility.  Therefore, our 
reducibilities have a more topological rather than 
computational character.  In particular, we define 
$f \leqT g $ to mean $f \leq_W^c \hat g$, 
and make similar topological definitions for 
$\leqtt$ and $\leqm$ in Section \ref{sec:Definition_of_tt_and_m_reducibility}.
We plan to address the 
question of the lightface theory in future work.

\section{Preliminaries}
\label{preliminaries}

\subsection{Notation}

We use standard computability-theoretic notation.  Brackets
$\langle m, n \rangle$ denote a canonical pairing function 
identifying $\omega\times\omega$ with $\omega$.  The 
expression $0^\omega$ refers to an $\omega$-length string 
of $0$'s.  Concatenation of finite or infinite strings $\sigma$
 and $\tau$ is denoted by 
 $\sigma\concat\tau$, which may be shortened to 
 $\sigma\tau$ in cases where it would cause no confusion. 
 If $\sigma$ is a string with a single entry $n$, 
 we also denote concatenation by $n\concat \sigma$ 
 or $\sigma \concat n$.  We usually use $X$ and $Y$ 
 to denote compact separable metric spaces, 
 $A,B,Z,W$ to denote elements of $2^\omega$ 
 or $\omega^\omega$, $C,D,P,Q$ to denote 
 subsets of $2^\omega$, and $\mathcal C, \mathcal D$ 
 to denote subsets of $\mathcal P(X)$. 
 Usually $f,g,$ and $j$ are arbitrary functions from $X$ to 
 $\R$ (the ones whose complexity we seek to categorize), 
 while $h,k,u,v,H$ and $K$ are typically continuous functions 
 from $\omega^\omega$ to $\omega^\omega$.  

\subsection{Computability and descriptive set theory}

We assume the reader is familiar with Kleene's $\mathcal O$ (but 
without it, one could still understand the results at the 
finite levels of each of the hierarchies).  The standard 
reference on this subject is Sacks \cite{sacksHRT}.
The $n$th jump 
of a set $A \in \omega^\omega$ is denoted $A^{(n)}$.  
For any $a \in \mathcal O^A$, 
if $|a| = n$ then let $A_{(a)}$ denote $A^{(n)}$, and if 
$|a|$ is infinite then let $A_{(a)}$ denote $H_{2^a}^A$.
If $a \in \mathcal O$ with $|a| = \alpha$, we will 
often simply write $\alpha$ instead of $a$.  Thus an expression 
like $A_{(\alpha)}$ is technically ambiguous, but since 
all the sets which it could refer to are one-equivalent, 
no problems will arise.

The reason for numbering the jumps in this 
lower-subscript way is to make 
them align correctly with the Borel hierarchy.  Recall 
that a set is $\mathbf \Sigma^0_1$ if it is open, 
$\mathbf \Pi^0_\alpha$ if it is the complement of a 
$\mathbf \Sigma^0_\alpha$ set, and $\mathbf \Sigma^0_{\alpha}$ 
if it is of the form $\cup_{n \in \omega} C_n$ 
where each $C_n$ is $\mathbf \Pi^0_{\beta_n}$ for some 
$\beta_n < \alpha$.  Then a set 
$C \subseteq \omega^\omega$ 
is $\mathbf \Sigma^0_\alpha$ if and only if there 
is a parameter $Z \in \omega^\omega$ and an index $i$ 
such that for all $A \in \omega^\omega$, 
$$A \in C \iff i \in (A \oplus Z)_{(\alpha)}$$
(and if no parameter is needed, we say $C$ 
is $\Sigma^0_\alpha$).

Still, at least once we will want to refer to the 
sets $H_a^A$, where $|a|$ is a limit ordinal.  In this case, 
we write $A^{(a)}$ to denote $H_a^A$.

\subsection{Representations}
\label{reps}

Although our results were motivated by considering $f:[0,1]\to {\mathbb R}$, 
they are also applicable in a wider context, \emph{represented spaces},
and hence, for completeness, we will briefly give an account of such spaces.
A standard reference is Weihrauch \cite{Weihrauch2000}.

In order for a machine to interact with a mathematical object, 
the object must be coded in a format a machine can read, 
such an element of $2^\omega$ or $\omega^\omega$.
For example, an element of $\R$ could be coded by 
a rapidly Cauchy sequence of rational numbers (which 
is itself coded by an element of $\omega^\omega$ using 
some fixed computable bijection $\omega \leftrightarrow \Q$).
It is not too hard to see that a similar 
method will also work for any computable metric space, where 
the role of the rationals is taken by (codes for) a computable dense subset.

A \emph{representation} of a space $X$ is a partial function $\delta:\subseteq \omega^\omega\to X$, so that 
elements $x\in X$ have $\delta$-names $A_x$ (strictly a set $\{A_x\mid \delta(A_x)=x\}$). Note that $x$ can have many names $A_x$, 
and not every element of $\omega^\omega$ is a name.  A representation 
induces a topology on $X$, the final topology, defined by $U \subseteq X$ 
is open if and only if $\delta^{-1}(U)$ is open in the subspace 
topology on $\dom \delta$.  If $X$ already has a topology, we 
restrict attention to representations which induce the topology of $X$.
Then if $X$ and $Y$ are represented spaces and 
$f:X\rightarrow Y$, we say $f$ is \emph{computable} if 
there is a computable function $F:\omega^\omega\rightarrow\omega^\omega$
such that whenever $A_x$ is a name for $x$, 
then $F(A_x)$ is a name for $f(x)$.  We say that $F$ \emph{realizes}
$f$.  Because $x$ and $f(x)$ each have many names, 
in general realizers are not unique.

 Not all representations are created equal.
 For example, the base 10 representation for reals is 
 a valid representation according to the above definition, but
the function $f(x)=3x$ is not computable with respect the 
base 10 representation on both sides (what digit
should the algorithm output first when seeing input .33333...?).  
However,
it is computable with respect to the Cauchy 
name representation on both sides.  This difference is 
captured in the following definition:  a representation 
$\delta:\subseteq \omega^\omega \rightarrow X$ is \emph{admissible}
if for every other continuous $\delta':\subseteq \omega^\omega\rightarrow X$, 
there is a continuous function $G:\omega^\omega\rightarrow\omega^\omega$
such that for all $A \in \dom \delta'$, we have $\delta(A) = \delta'(G(A))$.
That is, $G$ transforms $\delta'$-names to $\delta$-names.
Observe that it is possible to continuously 
transform a base 10 name for $x$ into a 
Cauchy name for $x$, but not vice versa.
Some definition chasing shows that 
the Cauchy name representation for $\R$ is admissible.
Restricting attention to admissible representations
allows continuity properties of $f$ to be reflected in 
its realizers.

\begin{thm}[Kreitz and Weihrauch \cite{KreitzWeihrauch1985}, Schr\"oder \cite{Schroeder2002}]
If $X$ and $Y$ are admissibly represented separable $T_0$ spaces, then 
a partial function $f:\subseteq X\to Y$ has a continuous realizer 
if and only if $f$ is continuous.\end{thm}

All of the pain and suffering 
involving representations is rewarded when we want to 
compare functions $f$ and $g$ in topologically incompatible areas, like 
Cantor space and ${\mathbb R}$. When comparing 
$f:X\to Y$ and $g:U\to V$, we can so so \emph{via their representations
in $\omega^\omega$.}
Given two represented
spaces $X$ and $Y$, a \emph{Weihrauch problem} is a
multivalued partial function $f: X \rightrightarrows Y$.
The $\rightrightarrows$ indicated that this definition 
concerns multivalued partial functions.
But 
in this paper, almost all problems will be total and
single-valued.  Accordingly, we will
freely call problems ``functions''.

We conclude this section with the precise definition
of Weihrauch reducibility on represented spaces.
Let $X,Y,Z$ and $W$ be represented spaces with
representations $\delta_X, \delta_Y, \delta_Z$,
$\delta_W$.
If $f: X \rightarrow Y$ and
$g: Z \rightarrow W$ are two
single-valued Weihrauch problems, we say
$f$ is \emph{Weihrauch reducible} to $g$, written $f \leq_W g$,
if there are computable functions 
$H,K:\subseteq\omega^\omega\rightarrow \omega^\omega$ such that
$f(\delta_X(A)) = \delta_Y H(A,B)$ for all $A \in \dom \delta_X$
and all $B$ such that $\delta_W B = g(K(A))$.
We say $f$ is \emph{strongly Weihrauch reducible}
to $g$, written $f \leq_{sW} g$, if $H(A,B)$ can be
replaced by $H(B)$ above.  The notions of
\emph{continuous Weihrauch reducibility} and
\emph{continuous strong Weihrauch reducibility},
denoted $\leq_{W}^c$ and $\leq_{sW}^c$ respectively,
are obtained by allowing $H$ and $K$ to be
merely continuous
rather than computable.
The parallelization $\hat g:Z^\omega\rightarrow W^\omega$
is defined by $\hat g( (z_i)_{i\in\omega}) = (g(z_i))_{i\in \omega}$.

In this paper, we will be dealing for the most part with 
situations where the coding is clear,
and hence suppress the $\delta_X$ notation whenever possible.

\subsection{Baire functions}

Baire functions are the most tractable functions we might consider after continuous ones.
Baire 1 functions are those which are defined as pointwise limits 
of a countable collection of continuous functions; $f(x)=\lim_s f_s(x)$ 
with each $f_s$ continuous.
More generally,
let $X$ be a compact separable metric space.
By $\mathcal C(X)$, we mean 
the continuous functions $f:X\rightarrow \mathbb R$.  The Baire 
hierarchy of functions on $X$ is defined as follows.  Let 
$\mathcal B_0(X) = \mathcal C(X)$.  For each $\alpha>0$, let 
$\mathcal B_{\alpha}(X)$ be the set of functions which are
pointwise limits of sequences of functions from 
$\cup_{\beta <\alpha} \mathcal B_\beta(X)$.  The functions in 
$\mathcal B_\alpha(X)$ are also referred to as the Baire $\alpha$ 
functions when $X$ is clear.

It is well-known that a function $f$ is Baire $\alpha$ if and 
only if the inverse image of each open set under $f$ is
$\mathbf\Sigma^0_{\alpha+1}$.  When $X =2^\omega$, the Baire $\alpha$ 
functions can also be characterized via the jump.

\begin{prop}[Folklore]\label{prop:baire-and-jump}
For each ordinal $\alpha$ and 
$f:2^\omega \rightarrow \mathbb R$, $f \in \mathcal B_\alpha(2^\omega)$ 
if and only if there is a Turing 
functional $\Gamma$ and $B \in 2^\omega$ such that 
$$f(A) = \Gamma((A\oplus B)_{(\alpha)}).\footnote{Technically 
$\Gamma$ outputs a code for $f(A)$ using some admissible representation.}$$
\end{prop}
\begin{proof}
When such $\Gamma$ and $B$ exist, one can 
readily check that the inverse images of 
open sets are $\mathbf \Sigma^0_{\alpha+1}$.  Conversely, 
if $f$ is Baire $\alpha$, then the sets $f^{-1}((p,q))$ 
for $p,q \in \mathbb Q$ can each be written as
$$f^{-1}((p,q)) = \{A : i_{p,q} \in (A \oplus B_{p,q})_{(\alpha+1)}\}$$
Therefore, if $B$ is an oracle containing each $B_{p,q}$ and $i_{p,q}$ 
in a uniformly accessible manner, one can use a 
$(A\oplus B)_{(\alpha)}$ oracle to enumerate the rational intervals 
$(p,q)$ containing $f(A)$, which is enough to make a Cauchy name 
for $f(A)$.
\end{proof}

\subsection{Ranks on Baire 1 functions}

In \cite{KechrisLouveau1990},
Kechris and Louveau defined
three ranks $\alpha$, $\beta$ and $\gamma$ on the Baire 1
functions.  These ranks had been used either explicitly or
implicitly in the literature analyzing this class of functions.
Given a Baire 1 function $f:X\rightarrow \mathbb R$, the
following derivation process is used to define the $\alpha$
rank.  Given rational numbers $p,q$ with $p<q$ and
a closed set $P \subseteq X$, let
$$P'_{p,q} = P \setminus \cup \{U \subseteq X : U \text{ is open and }
f(U) \subseteq (p, \infty) \text{ or } f(U) \subseteq (-\infty,q)\}$$

For a fixed pair $p,q$, define an $\omega_1$-length sequence
$\{P_\nu\}_{\nu < \omega_1}$ as follows.  Let $P_0 = X$,
$P_{\nu + 1} = (P_\nu)'_{p,q}$, and $P_{\nu} = \cap_{\mu < \nu} P_\mu$
if $\nu$ is a limit ordinal.  Since $X$ is separable, it has a
countable basis, so the sequence must stabilize below $\omega_1$.
Let $\alpha(f,p,q)$ be the least $\nu$ such that $P_\nu = \emptyset$;
one can show that such $\nu$ exists if and only if $f$ is Baire 1.

Finally, the $\alpha$ rank is defined by
$\alpha(f) = \sup_{p<q} \alpha(f,p,q)$.  The $\beta$ and $\gamma$
ranks are also defined by different
transfinite derivation processes.  Kechris and Louveau show that 
the levels of the form $\omega^\nu$ are especially robust
in the following sense.

\begin{thm}[\cite{KechrisLouveau1990}]
For any countable $\xi$ and any bounded 
Baire 1 function $f$, 
$$\alpha(f) \leq \omega^\xi \iff \beta(f) \leq \omega^\xi \iff
\gamma(f) \leq \omega^\xi.$$
\end{thm}

\section{Topological Turing reducibility on $2^\omega$}

First we define the topological Turing reducibility as mentioned in the
introduction.  First we fix $X = \C$.

\begin{defn}  For $f,g : \C \rightarrow \mathbb R$,
  let $f \leqT g$ if $f \leq_{W}^c g$.
\end{defn}

Equivalently, $f \leqT g$ if and only if there is a countable
sequence of continuous functions $k_i:2^\omega \rightarrow 2^\omega$
and a continuous function $h:\subseteq 2^\omega \rightarrow 2^\omega$
such that whenever $\{B_i\}_{i<\omega}$ are Cauchy names for
$\{g(k_i(A))\}_{i<\omega}$, $h(A \oplus \bigoplus_{i<\omega} B_i)$
is a Cauchy name for $f(A)$.  Observe that all continuous
functions are equivalent under $\leqT$.

The restriction of the domain to $\C$ is not essential, but helps 
keep the notation manageable.  If $X$ is a compact 
separable metrizable space and $f:X\rightarrow \mathbb R$, 
then in order to compare $f$ with other functions, we may
replace $f$ with $fd_X : 2^\omega \rightarrow \mathbb R$, 
where $d_X :2^\omega \rightarrow X$ is any admissible 
representation.  This gives a well-defined extension of the 
notion of $\leqT$ because, as the following proposition 
makes explicit, it does not matter which admissible 
representation we choose.

\begin{prop}\label{prop:admissible}
Let $X,Y$ be compact separable metrizable spaces and 
let $d_X, d_X':2^\omega\rightarrow X$ and 
$d_Y,d_Y':2^\omega\rightarrow Y$ be any admissible representations
for $X$ and $Y$ respectively.  Let $f:X\rightarrow \mathbb R$
and $g:Y\rightarrow \mathbb R$.  Then 
$$fd_X \leqT gd_Y \iff fd_X' \leqT gd_Y'$$
\end{prop}
\begin{proof}
Suppose that $(k_i)_{i<\omega}$ and $h$ witness that 
$fd_X \leqT gd_Y$.  By admissibility, let $\phi, \psi:2^\omega\rightarrow 2^\omega$ be continuous functions such that 
$d_X' = d_X \phi$ and $d_Y = d_Y' \psi$ (note the asymmetry).
Then $(\psi k_i \phi)_{i<\omega}$ and $h$ 
witness that $fd_X' \leqT gd_Y'$.
The reverse implication follows by symmetry.
\end{proof}

Therefore, from here on we may restrict our attention 
to functions $f,g :2^\omega \rightarrow \mathbb R$.

Note that one could also consider the notion defined by $f \leq_{sW}^c g$.
However, this is almost the same notion as the one defined.
If $f \leqT g$ via $\{k_i\}$ and $h$, and if
$g$ is a non-constant function, then letting $B_0$ and $B_1$
be such that $g(B_0) \neq g(B_1)$, one could additionally
consider the continuous functions $\{k_i'\}$ which map
$A$ to $B_0$ if $A(i) = 0$ and map $A$ to $B_1$ otherwise.
Then $A$ itself is continuously recoverable from
$\bigoplus_i k_i'(A)$, so by adding these to the original
$\{k_i\}$, a small modification to the original $h$
will do the job in the strong Weihrauch setting.

Therefore,
if $g$ is non-constant, then
$f \leq_{sW}^c g $ if and only $f \leqT g$.
On the other hand, if $g$ is constant, then
$\{f : f\leq_{sW}^c g\}$
is just the set of constant functions.  So there is
no need to consider the strong variant separately.

Now let us define some jump functions to characterize the
$\leqT$ degrees of the Baire functions.  The jump functions 
we consider are real-valued,
because of our original motivation to study functions
from $[0,1]$ to $\mathbb R$.  But the jump operator
can be represented as a real-valued function in a
standard way.

\begin{defn} For $n \in \omega$,
let $j_n:\C\rightarrow \R$ be defined by 
$$j_n(A) = \sum_{i \in A_{(n)}} 2^{-(i+1)}.$$
\end{defn}

Because each $j_n(A)$ is irrational, its binary expansion
can be continuously recovered from it.  Therefore,
by Proposition \ref{prop:baire-and-jump}, if $f$
is Baire $n$, then $f \leqT j_n$.  We can also extend
the definition to the ordinal notations.  Context will
make it clear whether natural number subscript
should be interpreted as a natural number or as an
ordinal notation.

\begin{defn} For $a \in \kO$, let 
let $j_a:\C\rightarrow \R$ be defined by 
$j_a(A) = \sum_{i \in A_{(a)}} 2^{-(i+1)}$.
If $a$ is a limit notation, let
$j^a:\C\rightarrow \R$ be defined by
$j^a(A) = \sum_{i \in A^{(a)}} 2^{-(i+1)}$.
\end{defn}

Therefore, if $f$ is Baire $\alpha$, then $f \leqT j_a$ for
all $a$ with $|a|=\alpha$.  The following properties are clear.

\begin{prop}\label{prop:jump_relations}
  For any notations $a,b \in \kO$,
  \begin{enumerate}
  \item $j_a \leqT j_b$ if and only if $|a| \leq |b|$.
  \item If $a$ and $b$ are limits with $|a|=|b|$, then $j^a \equiv_{\mathbf T} j^b$.
  \item If $a$ is a limit, $j^a <_{\mathbf T} j_a$.
  \end{enumerate}
\end{prop}
\begin{proof}
  All parts of the proposition which claim that a reduction exists
  follow from the fact that $|a|\leq |b|$ implies $H_a^A \leq_T H_b^A$,
  uniformly in $A$.  For the non-reductions, suppose for the sake 
  of contradiction that $j_a \leqT j_b$ with $|a|>|b|$ 
  or $j_a \leqT j^a$.  Let $Z\in2^\omega$ 
  be an oracle strong enough to compute 
  the continuous functions $\langle k_i\rangle$ and $h$ used 
  in the reduction.  Then $H_{2^a}^Z \leq_T H_{2^b}^Z$ 
  or $H_{2^a}^Z \leq_T H_a^Z$, which are not possible.
\end{proof}

The previous proposition justifies the use of notation $j_\alpha$ 
to refer to $j_a$ for some unspecified $a \in\kO$ with $|a| = \alpha$.
By relativization, we can go further up the ordinals.

\begin{defn} For any $Z \in 2^\omega$ any any $a \in \kO^Z$, 
define $$j_a^Z(A) =  \sum_{i \in (A \oplus Z)_{(a)}} 2^{-(i+1)},$$ 
and similarly for $j^{a,Z}$.
\end{defn}

Proposition \ref{prop:jump_relations} can then be generalized to  
replace $j_a$ and $j_b$ with $j_a^Z$ and $j_b^W$, 
under the assumption that 
$a,b \in \kO^Z \cap \kO^W$.
We leave both the statement and proof of this 
generalization to the reader, but 
for example, part (1) follows from the fact that 
$H_a^{(A \oplus Z)} \leq_T H_b^{(A\oplus Z) \oplus W}$ 
uniformly in $A$; in the generalization the forward 
reduction is the continuous map $A \mapsto A\oplus Z$, rather 
than the identity map as it was in the original.  Therefore, 
for any $\alpha < \omega_1$, we may use $j_\alpha$ 
to refer to $j_a^Z$ for some pair $Z,a$ with $Z \in 2^\omega$ 
and $a \in \kO^Z$ with $|a|^Z_{\kO} = \alpha$, and 
it does not matter which such $Z,a$ we use because they 
are all in the same $\leqT$ equivalence 
class.\footnote{Given $a,Z$ and $b,W$ with $|a|^Z_\kO = |b|^W_\kO$ 
but $a \not\in \kO^W$ or $b \not\in \kO^Z$, first fix $a' \in \kO^Z\cap \kO^W$ 
with $|a'|^Z_{\kO} = |a'|^W_\kO = |b|^W_\kO$, 
then observe $j_a^Z \leqT j_{a'}^Z \leqT j_{a'}^W \leqT j_b^W$ .}  Similar 
remarks apply to the expression $j^\alpha$.

We conclude by showing that every Baire function in $\mathcal F(\C,\R)$ 
is topologically Turing equivalent to one of the $j_\alpha$ or $j^\alpha$.
To reduce the notational clutter, we prove the version where 
$\alpha$ is constructive, and leave the relativization to the reader.

\begin{prop}\label{P2}
Let $\alpha$ be a constructive ordinal.  
If a Baire function $f$ is not Baire $\alpha$, then $j_{\alpha+1} \leqT f$.
If $\alpha$ is a limit and $f$ is not Baire $\beta$ for any $\beta < \alpha$,
then either $f  \equiv_{\mathbf T}  j^\alpha$, or 
$j_\alpha \leqT f$.
\end{prop}

\begin{proof}
Since $f$ is not Baire $\alpha$, there is 
an open set $U\subseteq \R$ such that $f^{-1}(U)$ 
is not ${\mathbf \Sigma}^0_{\alpha+1}$. 
Since $f$ is Baire, $f^{-1}(U)$ is Borel, 
so by Wadge determinacy $C_{\alpha+1} \leq_w f^{-1}(U)$,
where $C_{\alpha+1}$ is a canonical complete 
${\mathbf\Pi^0_{\alpha+1}}$ subset of 
$\C$, and $\leq_w$ is Wadge reducibility.  Let $v$ be a 
continuous function such that for all $Z$, 
$$Z \in C_{\alpha+1} \iff v(Z) \in f^{-1}(U).$$  We now show how to 
reduce $j_{\alpha+1}$ to $f$.  It suffices to be able to compute each 
bit of $A_{(\alpha+1)}$ on input $A$.  Given $A$ and $i$, 
uniformly compute $Z$ such that 
$i \not\in A_{(\alpha+1)}$ if and only if $Z \in C_{\alpha+1}$.  Expressing 
$i \in A_{(\alpha+1)}$ as the statement 
$\exists k [u (i, k ) \in A_{(\alpha)}]$ for some computable $u$, 
compute also a sequence 
$Z_k$ such that $u( i, k) \in A_{(\alpha)}$ if and only if 
$Z_k \in C_{\alpha+1}$.  Then asking for the values of
$f(v(Z))$ and $f(v(Z_k))$, wait until 
you see one of these enter $U$.  This proves the first part.

Now suppose that $\alpha$ is a limit,
$\alpha = \lim_n \alpha_n$.  If $f$ is not Baire $\beta$ for any
$\beta<\alpha$, then $f \geqT j_{\alpha_n}$ for each $n$.  From
this it is clear that $f \geqT j^\alpha$.  Suppose that there is
an open set $U$ such that $f^{-1}(U)$ is not $\mathbf \Sigma^0_\alpha$.
Then by the same argument as above, $f \geqT j_\alpha$.  On
the other hand, if $f^{-1}(U)$ is $\Sigma^0_\alpha$ for each open
$U$, then $j^\alpha \geqT f$ as follows.  Let $W$ be an oracle
such that $\{(A,p,q) : f(A) \in (p,q)\}$ is $\Sigma^0_\alpha(W)$.
Given access to the oracle $j^\alpha(A\oplus W)$, we can enumerate
$\{(p,q) : f(A)\in (p,q)\}$.  This suffices to compute $f(A)$.
\end{proof}

So that is the complete picture for $\leqT$.
The particularly strong way in which each Baire $\alpha$ 
function is reducible to $j_{\alpha}$ is in fact a continuous Weihrauch 
reduction.  However, the reduction of Proposition \ref{P2} is 
not a continuous Weihrauch reduction since we query different 
values of $f$ for each bit of $A_{(\alpha+1)}$.  So the parallelization is
certainly used.

Kihara has subsequently obtained a further characterization 
of the degree structure of $\leqT$ in terms of Martin reducibility on 
uniformly Turing order preserving operators; see \cite{kihara1}.

\section{Definition of topological $tt$- and $m$-reducibilities}\label{sec:Definition_of_tt_and_m_reducibility}

The classical notions of $tt$- and $m$-reducibility on infinite binary 
sequences operate by restricting the number of bits of the oracle 
used and the manner in which they are used.  In the case of a 
$tt$-reduction, in order to get the $n$th bit of the output, 
one specifies in advance, using only the number $n$, 
finitely many bits of the oracle that will be queried.  For each 
possible way the oracle could respond, one commits to 
an output for the $n$th bit.  Only then is the oracle queried and 
the commitment carried out.  The $m$-reducibility is even more 
restrictive.  In order to get the $n$th bit of the output, one specifies 
in advance a single bit of the oracle to query, and commits to 
copy the whatever the oracle has there as the $n$th bit.

As explained in the introduction, we have adopted the convention that
one bit of information about $f$ is an $\eps$-approximate answer to the
question ``Is $f(A)$ greater or less than $p$?'' where 
$A \in 2^\omega$ and $p \in \mathbb Q$.

Given $A \in 2^\omega$, $p \in \mathbb Q$,
and $\eps \in \mathbb Q^+$, we define the question 
$$\text{``}f(A) \lesssim_\eps p\text{''?}$$
so that ``yes'' or ``1'' is a correct answer if $f(A) < p+\eps$
and ``no'' or ``0''
is a correct answer if $f(A) > p-\eps$.  Observe that 
either answer is considered correct if $f(A)$ is within 
$\eps$ of $p$.  

We then define a representation of $\mathbb R$ 
whose domain is a subset of $2^\omega$, where each 
bit of a name for $y \in \mathbb R$ corresponds to 
a correct answer to a question of the form $y \lesssim_\eps p$.

\begin{defn}
We say $A \in 2^\omega$ is a \emph{separation name} for $y \in \mathbb R$ 
if for every $p\in \mathbb Q, \eps \in \mathbb Q^+$, we have 
$A(\langle p, \eps\rangle)$ correctly answers $y \lesssim_\eps p$.
\end{defn}

One can verify that the function mapping separation names to reals is 
an admissible representation.  
Now if we take the definition of $\leqT$ from the previous section, 
use the above representation for real numbers, and further specify that  
$h$ be either an $m$-reduction or a $tt$-reduction 
respectively, we obtain the following topological 
definitions of $\leqm$ and $\leqtt$.

\begin{defn}
We say $f \leqm g$ if and only if for every pair of rationals $p,\eps$, 
there are rationals $q,\delta$ and 
a continuous function $k:\C\rightarrow \C$ 
such that whenever $b$ is a correct answer 
to $g(h(A))\lesssim_\delta q$, 
$b$ is also a correct answer to $f(A)\lesssim_\eps p$.
\end{defn}

\begin{defn}
We say $f\leqtt g$ if and only if for every pair of rationals $p,\eps$, 
there are 
\begin{itemize}
\item finitely many rationals $(q_i, \delta_i)_{i<r}$
\item continuous functions $k_i:2^\omega\rightarrow 2^\omega$, and
\item a truth table function $h:2^r \rightarrow \{0,1\}$
\end{itemize}
such that whenever $\sigma \in 2^r$ is a string where each $\sigma(i)$ 
correctly answers $$g(k_i(A)) \lesssim_{\delta_i} q_i,$$ 
then $h(\sigma)$ correctly answers 
$f(A) \lesssim_\eps p$.
\end{defn}

It is clear that the reducibilities $\leqm$ and $\leqtt$
are reflexive and transitive, and that 
$$f\leqm g \implies f\leqtt g \implies f \leqT g.$$

Exactly as in 
Proposition \ref{prop:admissible}, these reductions 
may be more generally applied to functions whose domain is 
any compact separable metrizable space, using 
admissible representations.

Finally, all these reductions are primarily suitable for 
comparing discontinuous functions.

\begin{prop}\label{prop:4}
If $f$ is continuous and $g$ is non-constant, then $f \leqm g$.  
\end{prop}
\begin{proof} Since $g$ is non-constant, let $B_0, B_1 \in 2^\omega$ 
be such that $g(B_0)<g(B_1)$.
Given $p,\eps$, let $k$ be a continuous function which is 
equal to $B_0$ on $f^{-1}((-\infty, p-\eps])$ and 
equal to $B_1$ on $f^{-1}([p+\eps, \infty))$.  Since $f$ is continuous, 
these sets are closed, so such a $k$ exists.
Let $q, \delta$ be such that $g(B_0) < q-\delta < q + \delta < g(B_1)$. 
Then $q,\delta, k$ 
satisfies the part of the ${\mathbf m}$-reduction associated to $p,\eps$.
\end{proof}

\subsection{Equivalent definitions}\label{subsec:equivalent}
After hearing these results,
the following equivalent definitions for $\leqtt$ and $\leqm$ reducibilities 
were observed by Arno Pauly and Takayuki Kihara, respectively.

First some standard notation.  If $g:\subseteq X\rightrightarrows Y$ 
is a Weihrauch problem, $g^n$ is defined as the problem 
$g \times g: X^n \rightrightarrows Y^n$ where 
$(y_0,\dots y_{n-1}) \in g^n(x_0,\dots, x_{n-1})$ 
if and only if $g(x_i) = g(y_i)$ for all $i<n$.  
Then $g^\ast$ is defined as 
$g^\ast :\subseteq \cup_n X^n \rightrightarrows \cup_n Y^n$ 
where $\bar y \in g^\ast(\bar x)$ if $\bar x$ and $\bar y$ 
are the same length $n$ and $\bar y \in g^n(\bar x)$.

For any function $f:2^\omega\rightarrow \mathbb R$, let 
$S_f:\omega^\omega \rightrightarrows \{0,1\}$ be defined by
$$b \in S_f((p,\eps)\concat A) \iff b \text{ correctly answers } f(A) \lesssim_\eps p.$$

\begin{prop}[Pauly]
For $f,g:2^\omega\rightarrow \mathbb R$, $f \leqtt g$ if and only if 
$S_f \leq_{sW}^c S_g^\ast$.
\end{prop}
\begin{proof} If $g$ is constant, then each reducibility holds if and only if 
$f$ is constant as well.  So assume that $B_0,B_1 \in 2^\omega$ and 
$q,\delta \in \mathbb Q$ are 
inputs for which $g(B_0)<q-\eps < q+\eps < g(B_1)$.

If $f \leqtt g$, then for each $p,\eps$, 
let $(q_i,\delta_i,k_i)_{i<r}$ and $h$ be witness to this.  For each $p,\eps$, 
let $r'$ be the number of bits sufficient to describe $h$ according to some 
canonical self-delimiting coding.  Then define a strong Weihrauch reduction 
from $S_f$ to $S_g^\ast$ as follows:
\begin{itemize}
\item Given $(p,\eps)\concat A$, determine $r, r'$ from $(p,\eps)$ and set up 
a query to $S_g^{r'+r}$.
\item Use $r'$-many queries to ask about $(q,\delta)\concat B_0$ and 
$(q,\delta)\concat B_1$ in a sequence which encodes $h$.
\item Ask about $(q_i,\delta_i)\concat k_i(A)$ for each $i<r$.
\item Given the sequence of answers to these $r'+r$-many questions,
read off $h$ from the first $r'$ bits and apply it to the remaining $r$ bits.
\end{itemize}

The other direction uses the compactness of $2^\omega$.  Suppose that 
$S_f \leq_{sW}^c S_g^\ast$ via $K$ and $H$.  Fix $p$ and $\eps$.  
By compactness, there are finitely many strings $(\sigma_i)_{i<\ell}$ 
and for each $i$ there are finitely 
many rationals $(q_{ij},\delta_{ij})_{j<r_i}$ such that 
the cylinders $[\sigma_i]$ cover $2^\omega$, and 
for each $A \in 2^\omega$, if $\sigma_i \prec A$, then $K(A)$ 
has length $r_i$, and its $j$th coordinate 
begins with $(q_{ij},\delta_{ij})$.

Let $K_j$ be the function which computes the Cantor space part 
of the $j$th coordinate of $K$, when that coordinate exists.  That is, 
$K_j$ is defined by 
$$K((p,\eps)\concat \sigma_i C)(j) = (q_{ij},\delta_{ij})\concat K_j(\sigma_i C).$$
Let $(k_{ij})_{j<r_i}$ be functions 
do the following:
$$k_{ij}(A) = \begin{cases} B_0 & \text{ if } \sigma_i \not\prec A\\
K_j((p,\eps)\concat A) & \text{ if } \sigma_i \prec A. \end{cases}$$
Define also $k'_i(A) = B_j$ where $j=1$ if $\sigma_i \prec A$ and $0$ 
otherwise, and let $(q_i',\delta_i')$ be all equal to $(q,\delta)$.  
Let $r$ be the total number of $k_{ij}$ and $k'_i$ 
functions defined above.
Let $h: 2^{r}\rightarrow \{0,1\}$ 
be the truth table which uses the $k',q,\delta$ answers to determine 
which $\sigma_i \prec A$, then uses the $k_{ij},q_{ij}, \delta_{ij}$ answers to 
simulate the reverse reduction $H$.
\end{proof}


Kihara has also observed an equivalent definition of $\leqm$ related 
to partial order valued Wadge reducibility.  His definition and 
analysis also suggested a close variant of $\leqm$ whose theory 
may be even more natural than the one defined here.  We refer the reader to 
\cite{kihara1} for details.

\section{Properties of $\leqm$}\label{sec:properties_of_leqm}

In this section we prove our first main result concerning the 
$\leqm$ degrees of the jump functions $j_\alpha$ within the 
Baire $\alpha$ functions.
We start with some easier facts about 
the structure of the $\leqm$ degrees.  The proof of the following 
proposition is due to Kihara.

\begin{prop}\label{prop:game}
  For all Baire functions $f$ and $g$,
  we have either $f \leq_m g$ or $g \leq_m -f$.
\end{prop}
\begin{proof}
  We can understand the statement $f\leq_m g$ as saying that Player II
  has a winning strategy in the following game.  Player I 
  plays a target bit $\langle p,\eps \rangle$.  Player II plays
  its intended oracle bit $\langle q, \delta \rangle$.  Player I then
  starts playing bits of the input $A$; Player II also plays
  bits of a sequence $B$ in response, but Player II can pass (however
  they must ultimately produce an infinite sequence in order to win.)
  Player II wins if any correct answer to $g(B) \lesssim_\delta q$ 
  is also a correct answer to $f(A) \lesssim_\eps p$.  
  If Player II has a winning strategy, then
  $q, \delta$ and the continuous function
  $k$ defined by $k(A) = B$ are as in the definition of $\leqm$.
  But if Player I has a winning strategy, then for any $q, \delta$,
  there are $p, \eps$ (in fact, the same $p$ and $\eps$ each time,
  chosen according to the winning strategy of Player I) and a
  continuous function $k'$ which, following the winning strategy
  of Player I against Player II playing an arbitrary $B$,
  outputs $A = k'(B)$ such that either $g(B) < q + \delta$ and $f(A) \geq p + \eps$
  or $g(B) > q -\delta$ and $f(A) \leq p - \eps$.  Therefore, if $-f(A) < -p + \eps$,
  we must be in the first case and thus $g(B) < q + \delta$.
  Similarly, if $-f(A) > -p-\eps$ then we must be in the second 
  case, so $g(B) > q- \delta$.  This shows that 
  $g \leq_m -f$ via $k'$ (observe that $(-p, \eps)$ is the bit of $f(A)$
  actually queried).
\end{proof}

\begin{cor}
If $f \in \B_\alpha$, then $f \leqm j_{\alpha+1}$. 
\end{cor}
\begin{proof}
  If not, then by Proposition \ref{prop:game}
  we would have $j_{\alpha+1} \leqm -f \leqT j_\alpha$, impossible
  as $j_{\alpha+1}$ is not $\B_\alpha$.
\end{proof}

Our first theorem shows that 
the jump functions are the weakest functions in each Baire class.

\begin{thm}\label{T1} If $f$ is Borel
and $f \not\in \B_\alpha$,
then either $j_{\alpha+1} \leqm f$ 
or $-j_{\alpha+1} \leqm f$.
\end{thm}

It is easy to see why this theorem is true when $\alpha=0$.  If $f$ is 
not continuous, let $(z_n)_{n\in \omega}\rightarrow z$ be a convergent sequence 
of inputs for which $f(z) \neq \lim_n f(z_n)$. 
Without loss 
of generality, there is some $\delta>0$ such that for all $n$, $f(z_n) > f(z) + \delta$, 
or for all $n$, $f(z_n) < f(z) - \delta$.  In the first case, we have that 
$j_1 \leqm f$ via the following algorithm.  On input $(p,\eps)$, choose 
$(q,\delta')$ so that $[q-\delta', q+\delta'] \subseteq (f(z), f(z)+\delta)$.  
Then let $h$ be the function which, on input $x$, outputs bits of 
$z$ while computing approximations to $j_1(x)$.  If it ever sees 
that $j_1(x) > p-\eps$, it switches to outputting bits of $z_n$ 
for some $n$ large enough that $z$ and $z_n$ agree on all 
bits which were already committed to.    The case where 
$f(z_n) < f(z) - \delta$ for all $n$ is similar, only in that case 
we find that $-j_1 \leqm f$.

To prove this theorem in the general case 
we will make use of the following generalization of Borel Wadge 
determinacy.  We provide a simple proof of this 
generalization using Borel determinacy, but it is interesting 
to note that Louveau and Saint-Raymond \cite{LouveauSaintRaymond1987, LouveauSaintRaymond1988},
showed that this generalization is provable in
second order arithmetic via a much more intricate 
argument.  Therefore, 
the use of Borel determinacy here can be avoided.

\begin{prop}
  \label{thm:dichotomy} 
Let $D,E_0, E_1 \subseteq \omega^\omega$ be Borel. Then one of the following holds:
  \begin{enumerate}
  \item There is a continuous function $\varphi:\Baire \rightarrow \Baire$ such that $\varphi(D) \subseteq E_0$ and $\varphi(\Baire \setminus D) \subseteq E_1$.
  \item There is a continuous function $\psi:\Baire \rightarrow \Baire$ such that
    $\psi(E_0) \subseteq \Baire \setminus D$ and $\psi(E_1) \subseteq D$.
  \end{enumerate}
\end{prop}
\begin{proof}
  Define a two player game, where at turn $n$ player I (who plays first) plays $x(n)$ and player II plays $y(n)$. At the end of the game, II wins if
  \[(x \in D \wedge y \in E_0) \vee (x \not \in D \wedge y \in E_1).\]
By Borel determinacy, one of the two players has a winning strategy.  A winning strategy for II gives a continuous function meeting outcome (1).

  If on the other hand I has a winning strategy, then for every play of the game according to I's winning strategy we have that
  \[(x \not \in D \vee y \not \in E_0) \wedge (x  \in D \vee y \not \in E_1).\]
  This gives a continuous function meeting outcome (2).
\end{proof}

We give a new corollary to this theorem.

\begin{cor}
  \label{cor: reduction}
  Let $V \subseteq \Baire$ be $\mathbf \Pi^0_{\alpha}$.  Let $W \subseteq \Baire$ be  $\mathbf \Pi^0_{\alpha}$-hard and let $\{W_i\}_{i \in \Nat}$ be a partition of $W$ into Borel sets. Then there is a continuous function $\varphi\colon \Baire \rightarrow \Baire$ and $i \in \Nat$  such that:
  \begin{enumerate}
  \item $\varphi(V) \subseteq W_i$.
  \item $\varphi(\Baire \setminus V) \subseteq \Baire \setminus W$.
  \end{enumerate}
\end{cor}
\begin{proof}
  For each $i$, we can apply Theorem~\ref{thm:dichotomy} with $D=V$, $E_0 = W_i$ and $E_1 = \Baire \setminus W$. Assume that for each $i$, the second option of the theorem holds, i.e.\ there is a 
  continuous function $\psi_i$ such that 
  \[\psi_i(W_i) \subseteq \Baire \setminus V\mbox{ and }\psi_i(\Baire \setminus W) \subseteq V.\]
  Now take $K_i = \psi_i^{-1}(\Baire\setminus V)$. Note that $K_i$ is $\BS{\alpha}$ and  we have that $W = \bigcup_i W_i \subseteq \bigcup_i K_i$. Further, for all $i$ we know that $K_i \cap (\Baire \setminus W) = \emptyset$. Hence $W = \bigcup_i K_i$ and so $W$ is $\BS{\alpha}$. This is a contradiction as we are given  that $W$ is $\BP{\alpha}$-hard.

  Hence for some $i$ we have that the first option of Theorem~\ref{thm:dichotomy} holds.  That is, there is a
  continuous $\varphi$ such that $\varphi(V) \subseteq W_i$ and $\varphi(\Baire \setminus V) \subseteq \Baire \setminus W$.  
\end{proof}

\begin{proof}[Proof of Theorem~\ref{T1}]
  Since $f$ is not Baire $\alpha$, let $U\subseteq \mathbb R$ be an open set such that
  $f^{-1}(U)$ is not $\BS{\alpha+1}$. Without loss of generality, $U$ is of the form $(u, +\infty)$ or $(-\infty, u)$.
  If $U$ is of the form $(-\infty,u)$, 
then we will have $j_{\alpha+1} \leqm f$, and in the other case 
$j_{\alpha+1} \leqm -f$ (or equivalently, $-j_{\alpha+1} \leqm f$).  
Replacing $f$ with $-f$ if necessary 
let us assume $U = (-\infty, u)$.

Denote $f^{-1}(U)$ by $W$. 
Since $f$ is Borel,
$W$ is Wadge determined, so it is  
$\mathbf \Pi^0_{\alpha+1}$-hard. We can partition $W$ into the following sets
$W_0 = f^{-1}((-\infty, u-1])$ and for all $i \ge 1$,
$$W_i = f^{-1}\left(\left(u-\frac{1}{i}, u-\frac{1}{i+1}\right]\right).$$

  Take any $p, \epsilon \in \Q$ with $\epsilon >0$.
  Let $V = j_{\alpha+1}^{-1}((-\infty, p -\epsilon])$. The set $V$ is $\Pi^0_{\alpha+1}$. 
  (We have $A \in V$ if and only if for all finite $F\subseteq \omega$ such 
  that $\sum_{i\in F} 2^{-(i+1)} > p-\epsilon$, there is some $i \in F$ such that 
  $i \not \in A_{(\alpha+1)}$.  Recall from the introduction that $\{A : i \in A_{(\eta)}\}$ is 
  a $\Sigma^0_\eta$ set.)
  
    Thus by Corollary~\ref{cor: reduction} there is a continuous map $\varphi$ and an $i \in \Nat$ such that $\varphi(V) \subseteq W_i$ and $\varphi(2^\omega \setminus V) \subseteq f^{-1}([u, +\infty))$. Hence taking $\delta = \frac{1}{2(i+1)}$ and $q = u-\delta$ we have that for any $A \in 2^\omega$, there is only one correct answer to 
$f(\varphi(A)) \lesssim_\delta q$. Further, this is  also a correct answer to $j_{\alpha+1}(A)  \lesssim_\eps p$.
\end{proof}

\begin{cor}
If $g$ is Baire, $g\not\in \B_\alpha$ and $f \in \B_\alpha$, then $f \leqm g$.
\end{cor}

\section{The Bourgain rank on $\B_1$}

The structure of the $\leqm$-degrees and $\leqtt$-degrees 
within the Baire 1 functions 
is related to the $\alpha$ rank, also known as the Bourgain rank, 
which was studied by Kechris and 
Louveau \cite{KechrisLouveau1990}.  
Here we place that rank in a slightly more general setting 
that will be suitable for describing both the $\leqm$ and 
$\leqtt$ degrees, and establish some notation that will be 
used throughout.

\begin{defn}\label{D10a}
For any collection $\mathcal C \subseteq \mathcal P(X)$, 
a \emph{derivation sequence} 
for $\mathcal C$ is defined for $\nu<\omega_1$ by
\begin{itemize}
\item $P^0 = X$.
\item $P^{\nu+1} \supseteq P^\nu \setminus \cup\{U \text{ open } : \text{ for some }
C \in \mathcal C, P^\nu \cap U \subseteq C\}$
\item $P^\lambda \supseteq \cap_{\nu<\lambda} P^\nu$.
\end{itemize}
By replacing $\supseteq$ with $=$ in two places, we obtain the 
definition for the \emph{optimal derivation sequence} for $\mathcal C$.
\end{defn}

Here are some properties of derivation sequences which 
will be useful and which follow directly from the definitions.
\begin{prop}\label{prop:7} 
Let $Q^\nu$ be a derivation sequence for $\mathcal C
\subseteq \mathcal P(X)$.
\begin{enumerate}
\item If 
$P^\nu$ is the optimal derivation sequence for $\mathcal C$, then 
$P^\nu \subseteq Q^\nu$ for all $\nu$.  
\item If $k:X\rightarrow X$ is continuous, then $R^\nu:=k^{-1}(Q^\nu)$
is a derivation sequence for $\{k^{-1}(C) : C \in \mathcal C\}$.
\item If $\mathcal D\subseteq \mathcal P(X)$ 
is such that for every $C \in \mathcal C$, 
there is a $D \in \mathcal D$ such that $C \subseteq D$, 
then $Q^\nu$ is a derivation sequence for $\mathcal D$.
\end{enumerate}
\end{prop}

\begin{defn}[Bourgain rank, also known as $\alpha$ rank]\label{D10}
For $f \in \B_1$ and rationals $p,\eps$, let $P^\nu_{f,p,\eps}$ be the optimal 
derivation sequence for 
$\{ f^{-1}((-\infty, p+\eps)), f^{-1}((p-\eps,\infty))\}.$
Let $\alpha(f,p,\eps)$ be least ordinal $\nu$ such that 
$P^\nu_{f,p,\eps} = \emptyset$.
Let the Bourgain rank of $f$ be $$|f|_\alpha = \sup_{p,\eps \in \mathbb Q} \alpha(f,p,\eps).$$
\end{defn}

If $f, p, \eps$ are clear from context, we may write $P^\nu$ or $P^\nu_f$
instead of $P^\nu_{f,p,\eps}$. Observe that the compactness of $X$ implies 
that $\alpha(f,p,\eps)$ is always a successor, but in general $|f|_\alpha$ 
may be either a limit or a successor.

In the course of the optimal derivation process, individual points 
leave at various stages, and we would like to keep track of this.

\begin{defn} Let $A \in X$. If $P^\nu$ is the optimal derivation sequence for sets 
$\mathcal C$ and $P^\nu$ is eventually empty, let 
$|A|_\mathcal C$ denote the least $\nu$ such that 
$A \not\in P^\nu$.  
Given $f \in \B_1$, and $p,\eps$, let $|A|_{f,p,\eps}$ be the least 
$\nu$ such that $A \not\in P^\nu_{f,p,\eps}$.
\end{defn}

If $f,p,\eps$ and/or $\mathcal C$ are clear from context, we may 
just write $|A|_f$ or $|A|$.  Observe 
that $|A|$ is always a successor ordinal.  

The Bourgain hierarchy can be understood as a higher type 
version of the Ershov hierarchy.  Recall the Ershov hierarchy 
stratifies the $\Delta^0_2$ subsets of $\omega$ according 
to the amount of mind-changes needed in an optimal 
limit approximation to that set.  In general, ordinal-many 
mind-changes can be needed.  For $a \in \kO$,
a function $A: \omega\rightarrow \omega$ 
is $a$-computably approximable if there is a partial 
computable $\varphi(n,b)$ such that $A(n) = \varphi(n,b_n)$, 
where $b_n$ is the $\leq_\kO$-least ordinal 
$b_n \leq_\kO a$ for which the computation converges.  
We picture this process dynamically -- a computable 
procedure makes a guess about $A(n)$ associated to 
a certain ordinal.  If it changes its guess, it must decrease 
the ordinal.  This limits the number of mind-changes.

We can understand each open set removed as a part of the 
Bourgain derivation process as a guess about the 
answer to the question $f(x) \lesssim_\eps p$.  The 
open sets removed later in the derivation process,
have a high associated ordinal rank and correspond to 
early guesses; the open sets removed at the beginning 
of the derivation process correspond to the latest guesses.
The following object, a mind-change sequence, is nothing 
more than a derivation sequence annotated with the 
guesses that justified the derivation.  It can also be 
viewed as a higher-type analog of
$\varphi$ as above.
To simplify the notation, we assume
$\mathcal C = \{C_i : i < k\}$, where $k$ could be finite or $\omega$,
and $X = 2^\omega$.  Let $\Ord$ denote the ordinals.

\begin{defn}
Given $\mathcal C = \{C_i : i<k\} \subseteq \mathcal P(2^\omega)$, 
a \emph{mind-change sequence} for $\mathcal C$ is a
countable subset of
$M\subseteq \Ord \times 2^{<\omega} \times k$ for which
\begin{enumerate}
\item The sequence $Q^\nu$ defined by 
$$Q^\nu = 2^\omega \setminus \left( \bigcup_{\substack{(\mu,\tau,j)\in M\\ \mu < \nu}} [\tau] \right)$$
is a derivation sequence for $\mathcal C$, and
\item For all $(\nu, \sigma, i) \in M$, $[\sigma] \cap Q^\nu \subseteq C_i$.
\end{enumerate}
An \emph{optimal mind-change sequence} for $\mathcal C$ is 
one in which $Q^\nu$ is the optimal derivation sequence for $\mathcal C$.
\end{defn}

Observe that an optimal mind-change sequence always exists, 
since it just keeps track of the open sets $[\sigma]$ which are removed 
at stage $\nu$ of the construction of the optimal derivation 
sequence, and keeps track of which set $C\in\mathcal C$ caused 
$[\sigma]$ to be removed at stage $\nu$.

Two ``mind-change'' based encodings of the Baire 1 functions are 
suggested by this idea.  One encoding of $f \in \mathcal B_1$, 
following the $\alpha$ rank, would consist of a countable collection 
of mind-change sequences $M_{p,\eps}$, 
one for each $$\mathcal C_{p,\eps} = \{f^{-1}((p-\eps,\infty)), 
f^{-1}((-\infty, p+\eps))\}$$ for $p,\eps \in \mathbb Q$.
Another encoding, following the $\beta$ rank, would consist of
a different countable collection of mind-change sequences $M_\eps$, 
one for each
$$\mathcal C_\eps = \{f^{-1}((q-\eps,q+\eps)) : q \in \mathbb Q\}$$
for each $\eps\in\mathbb Q^+$.  We will not need to use such encodings 
explicitly, so we avoid further technical definitions, but this way 
of thinking about a Baire 1 function motivates all
the arguments which follow.

A mind-change sequence can serve as evidence of an upper
bound on the length of an optimal derivation sequence for
a collection $\mathcal C$.  The next notion provides evidence 
of a lower bound.  The idea is that if $[\sigma]$ was not 
removed at stage $\nu$ of the derivation process, then 
for each $C \in \mathcal C$, there was some element 
$A_{\nu,\sigma,C}$ which witnesses that 
$P^{\nu}\cap [\sigma] \not\subseteq C$.  
If $\mathcal C$ is countable 
and the derivation process lasts only countably 
many stages, then only countably 
many $A$ are needed to witness the necessity 
of an optimal derivation sequence being as long as it is.
Below, we define a \emph{scaffolding sequence} to be 
any countable collection of $A$'s which can supply 
all necessary witnesses, together with a record of where
in the process these $A$ are slowing things down.

\begin{defn}
Given $\mathcal C = \{C_i : i<k\} \subseteq \mathcal P(2^\omega)$, 
let $P^\nu$ be its optimal derivation sequence.  A 
\emph{scaffolding sequence} for $\mathcal C$ is any enumeration 
of a countable subset 
$S \subseteq 2^\omega \times \Ord \times 2^{<\omega} \times k$
such that
\begin{enumerate}
\item If $(A, \nu, \sigma, i) \in S$, then 
$A \in P^\nu \cap [\sigma] \setminus C_i$, and
\item If $P^\nu \cap [\sigma] \not\subseteq C_i$, there is $A \in 2^\omega$ 
with $(A, \nu, \sigma, i) \in S$.
\end{enumerate}
\end{defn}

Letting $S'$ be the projection of $S$ onto its first coordinate,
observe that for all $\mu<\nu$ and $\sigma$, if 
$P^{\nu} \cap [\sigma] \neq \emptyset$,
then $P^\mu \cap [\sigma] \cap S' \neq \emptyset$.


\section{Characterization of the $\leqm$ equivalence classes in $\B_1$}

In this section we prove parts (2)-(4) of Theorem \ref{thm:T2},
characterizing the structure of the $\leqm$ degrees within the Baire 1 
functions.

We will need to consider the case when $|f|_\alpha$ is a successor 
with special care. Supposing we have such an $f$, 
let $\nu, p,\eps$ be defined so that $\nu + 1 = \alpha(f,p,\eps) = 
|f|_\alpha$.  
Of course, we may also have $\nu + 1 = \alpha(f,p',\eps')$ for some 
other rationals $p',\eps'$.

\begin{defn}
Given $f \in \B_1$ with $|f|_\alpha=\nu+1$, and $p,\eps \in \mathbb Q$, 
say $(p,\eps)$ is \emph{maximal} if $f(P^\nu_{f,p,\eps}) \setminus (p-\eps,p+\eps) \neq \emptyset$
and $\alpha(f,p,\eps) = \nu+1$.  
\end{defn}

Observe that maximal $(p,\eps)$ always exist.  If 
$P^\nu_{f,p,\eps} \neq \emptyset$,
but $f(P^\nu_{f,p,\eps}) \setminus (p-\eps,p+\eps) = \emptyset$, then by 
decreasing $\eps$, one may shrink $(p-\eps,p+\eps)$ to include 
an element of $f(P^\nu_{f,p,\eps})$ (which grows in size).

\begin{defn}Let $f \in\mathcal B_1$ with $|f|_\alpha = \nu+1$. We say $f$ is
\begin{itemize}
\item \emph{two-sided} if there is a maximal $(p,\eps)$ 
such that $f(P^\nu) \not\subseteq (p-\eps,\infty)$ and 
$f(P^\nu) \not\subseteq (-\infty, p+\eps)$;
\item \emph{one-sided} otherwise;
\item \emph{left-sided} if for every maximal $(p,\eps)$,
$f(P^\nu) \subseteq (-\infty, p+\eps)$;
\item \emph{right-sided} if for every maximal $(p,\eps)$, 
$f(P^\nu) \subseteq (p-\eps,\infty).$
\end{itemize}
\end{defn}

For example $j_1$ is left-sided, as is any 
discontinuous lower semi-continuous 
function.  If $f$ is left-sided, then $-f$ is right-sided, 
and vice versa.  However, there are one-sided $f$ which are 
neither right-sided nor left-sided.  For example, consider 
$$f(X) = \begin{cases} 1 & \text{ if } X\in [0] \setminus \{01^\omega\} \\
-1 & \text{ if } X \in [1] \setminus \{10^\omega\} \\
0 & \text{ otherwise}.
\end{cases}$$

We are now ready to prove the parts (2)-(4) of our second main theorem.

\begin{thm}\label{T2}
For $f,g \in B_1$, $|f|_\alpha < |g|_\alpha$ implies $f\leqm g$. If 
$|f|_\alpha = |g|_\alpha$, then $f\leqm g$ if and only if 
at least one of the following holds:
\begin{enumerate}
\item $|f|_\alpha$ is a limit ordinal.
\item $g$ is two-sided.
\item $f$ is one-sided and $g$ is neither right-sided nor left-sided.
\item $f$ and $g$ are either both right-sided or both left-sided.
\end{enumerate}
\end{thm}
\begin{proof}

We begin with a general observation.  Suppose that $p,\eps, q, \delta \in 
\mathbb Q$ and $k:2^\omega\rightarrow 2^\omega$ is a continuous 
function such that any correct answer to $g(k(A)) \lesssim_\delta q$
is also a correct answer to $f(A) \lesssim_\eps p$. Then
for any $A$, $g(k(A)) < q+ \delta$ implies $f(A) < p + \eps$, 
so $k^{-1}(g^{-1}((-\infty, q+\delta)) \subseteq f^{-1}((-\infty, p + \eps))$.
Similarly, $k^{-1}(g^{-1}((q-\delta,\infty)) \subseteq f^{-1}((p-\eps,\infty))$.
Therefore, the sets $Q^\mu$ defined by 
$$Q^\mu = k^{-1}(P^\mu_{g,q,\delta})$$
are a derivation sequence for $\{f^{-1}((-\infty,p+\eps)), 
f^{-1}((p-\eps,\infty))\}$.  Therefore 
$P^\mu_{f,p,\eps} \subseteq Q^\mu$ for each $\mu$,
so $\alpha(f,p,\eps) \leq \alpha(g,q,\delta)$. 
Furthermore, for all $A \in 2^\omega$, we have 
$|A|_{f,p,\eps} \leq |k(A)|_{g,q,\delta}$.

Now suppose $f\leqm g$.  Then for any $p,\eps$, there are $q,\delta$
and $k$ as above, so $|f|_\alpha \leq |g|_\alpha$.  The first 
statement of the theorem now follows by Proposition \ref{prop:game}
and the observation that $|g|_\alpha = |-g|_\alpha$ for all $g$.

From now on we consider the case where $|f|_\alpha = |g|_\alpha$. 

Suppose that $f\leqm g$.  We claim that if 
$|f|_\alpha = \nu +1$ (a successor) then one of (2)-(4) 
in the statement of the theorem holds.

Let $(p,\eps)$ be maximal for $f$.  Since $f\leqm g$,
let $q,\delta$ and $k$ be as in the first paragraph.
By the choice of $p$ and $\eps$, there is an $A \in 2^\omega$ 
with $|A|_{f,p,\eps} = |k(A)|_{g,q,\delta} = \nu +1$
and $f(A) < p-\eps$, or there is a $B \in 2^\omega$ with 
$|B|_{f,p,\eps} = |k(B)|_{g,q,\delta} = \nu + 1$ 
and $f(B) > p +\eps$, or perhaps both occur.
If such $A$ exists, then $g(k(A)) < q-\delta$ and if such 
$B$ exists, then $g(k(B)) > q+\delta$.  

Therefore, 
if $g$ is not two-sided, then $f$ is not two-sided; in 
that case, if $g$ is right-sided or left-sided, then $f$ 
must match.  This completes the proof that $f \leqm g$ 
implies the disjunction of (1)-(4).

Assuming now the disjunction of (1)-(4),  
let $p,\eps$ be given.  First we choose a pair $q,\delta$ 
which gives us enough room to work.  
If $|f|_\alpha = \nu + 1$,
choose $(q,\delta)$ to be
maximal for $g$.  Additionally, if $g$ is two-sided, 
make sure $q$ and $\delta$ 
witness the two-sidedness of $g$.
Or, if $f$ is one-sided and $g$ is neither left-sided nor 
right-sided, then if 
$f(P^\nu_{f,p,\eps}) \subseteq (p-\eps, \infty)$ 
(respectively $f(P^\nu_{f,p,\eps}) \subseteq (-\infty, p+\eps)$)
make sure $g(P^\nu_{g,q,\delta}) \subseteq (q-\delta, \infty)$
(respectively $g(P^\nu_{g,q,\delta}) \subseteq (-\infty, q + \delta)$).
If $f$ and $g$ are both right- or both left-sided, 
a maximal choice of $q$ and $\delta$ suffices without 
further restrictions.
If $|f|_\alpha$ is 
a limit, choose $q,\delta$ so that 
$\alpha(f,p,\eps) < \alpha(g,q,\delta)$ and
$g(P^\nu_{g,q,\delta})\setminus (q-\delta,q+\delta) \neq \emptyset$ 
(decreasing $\delta$ if necessary to achieve the latter).
In this case, 
define $\nu$ so that $\alpha(g, q, \delta) = \nu + 1$.

We now define a continuous function $k$ such that 
any correct answer to $g(k(A)) \lesssim_\delta q$ 
also correctly answers $f(A) \lesssim_\eps p$.
Given $A$, its image $k(A)$ will be defined in 
stages according to an algorithm 
which uses oracle information about a mind-change sequence 
related to $f$ and a scaffolding sequence related to $g$.  
By defining $k(A)$ in stages, we guarantee $k$ is continuous.

Let $\mathcal C = \{C_0,C_1\}$, where $C_0 = f^{-1}((-\infty, p+\eps))$ 
and $C_1 = f^{-1}((p-\eps, \infty))$.  Let $\mathcal D = \{D_0,D_1\}$, 
where $D_0 = g^{-1}((-\infty, q+\delta))$ and 
$D_1 = g^{-1}((q-\delta, \infty))$.
Let $Z$ be an oracle which contains the following information:
\begin{itemize}
\item A well-order $W$ long enough that $\nu$ has a code in $\kO^W$
(a technical point which allows us to use $\kO^W$ in place of $\Ord$ in the 
mind-change and scaffolding sequences).
\item An optimal mind-change sequence $M$ for $\mathcal C$.
\item A scaffolding sequence $S$ for $\mathcal D$.
\end{itemize}

Letting $A$ denote the input, 
at each stage $s$, we will have defined an initial segment 
$\tau_s$ of $k(A)$.  We will be 
keeping track of an ordinal $\mu_s$, an index $i_s \in \{0,1\}$, 
and an element $B_s \in 2^\omega$, where $\tau_s \prec B_s$.  
We will always maintain the following:
\begin{enumerate}[(i)]
\item that $|A|_\mathcal C \leq \mu_s+1 \leq |B_s|_\mathcal D$, 
\item that $i_s$ is the only correct 
answer to $g(B_s)\lesssim_\delta q$, and
\item if $|A|_\mathcal C = \mu_s+1$, then $i_s$ correctly answers 
$f(A) \lesssim_\eps p$.
\end{enumerate}

The idea always is that as long as it seems like 
 $|A|_\mathcal C = \mu_s+1$, 
we are working towards making $k(A) = B_s$.   
If we later see the bound on $|A|_\mathcal C$ drop, 
and $i_s$ no longer looks like a suitable answer,
then because $|B_s|_\mathcal D$ is large, 
no matter how much of $B_s$ has been copied, we 
can switch to a nearby $B_t$ for which $i_t = 1-i_s$ is 
the only correct answer to $g(B_t) \lesssim_\delta q$, 
and $|B_t|_\mathcal D$ is still large.

Let 
$\lambda$ denote the empty string. 
Let $\tau_0 = \lambda$. 
We begin differently depending on whether 
$g$ is two-sided. In both of the following cases, 
the reader can verify that conditions (i)-(iii) 
are satisfied at stage $s=0$.

If $g$ is two-sided, we first 
wait until we see $A$ leave $P^{\mu}_{f,p,\eps}$ for 
some $\mu \leq \nu$.  That is, we see $(\mu, \sigma, i)$ in $M$ 
with $\sigma \prec A$.  Let $\mu_0=\mu$ and $i_0=i$.
Now, since $g$ is two-sided, 
regardless of $i$, $P^\nu_{g,q,\delta}\setminus D_{1-i}$ is non-empty,
and we can find an element $B$ in this set (by looking in $S$ for 
something of the form $(B,\nu,\lambda, 1-i)$). Let $B_0 = B$.

If $g$ is not two-sided, then for some $j$, 
$P^\nu_{g,q,\delta}\setminus D_j$ is non-empty, so first 
we wait until we see an element $B$ and a $j$ to witness this (by 
looking in $S$ for something of the form $(B,\nu,\lambda, j)$).
Let $\mu_0 = \nu$, $i_0 = 1-j$, and $B_0 = B$. By the choice 
of $q$ and $\delta$, if $|A|_\mathcal C = \nu +1$, then 
$i_0$ correctly answers $f(A)\lesssim_\eps p$.\footnote{
In case (1), 
by the choice of $\nu$, $P^\nu_{f,p,\eps} = \emptyset$, so 
$|A|_\mathcal C < \nu+1$.  In case (3), 
$f$ is one-sided, 
so $P^\nu_{f,p,\eps} \subseteq C_i$ for some $i$.  Note that 
in this case, we have chosen $q,\delta$ specifically to 
make sure that $j = 1-i$.  In case (4), we also have 
$P^\nu_{f,p,\eps} \subseteq C_{1-j}$ (note that $j=1$ if 
$f$ and $g$ are both left-sided and $j=0$ if $f$ and $g$ 
are both right-sided).}

At stage $s+1$, set $\mu_{s+1}$ to be the least $\mu$ for which 
we have seen $A$ leave $P^\mu_{f,p,\eps}$.  If $\mu_{s+1}<\mu_s$, 
that is because $(\mu_{s+1},\sigma, i)$ just entered $M$ 
for some $\sigma \prec A$.  If $i = i_s$, let $B_{s+1} = B_s$ 
and $i_{s+1} = i_s$. 
But if $i \neq i_s$, then set $i_{s+1} = i$, and look through 
$S$ to find a $B$ so that 
$$B \in P^{\mu_{s+1}}_{g,q,\delta}\cap[\tau_s]\setminus D_{i_s}.$$ 
Such a $B$ must exist because $B_s$ witnesses that 
$P^{\mu_s}_{g,q,\delta} \cap [\tau_s]$ is non-empty. Let 
$B_{s+1} = B$.  Finally, let $\tau_{s+1} = B_{s+1}\restrict |\tau_s|+1$. 
That completes the construction.

At each stage the properties (i)-(iii) are maintained.  
Now if $|A|_\mathcal C = \mu+1$, there is a stage $s$ 
at which it is seen that $A$ leaves $P^\mu_{f,p,\eps}$. 
The $\mu_s, i_s$ and $B_s$ defined at that stage never 
change again.  Then $k(A) = B_s$, and the only correct answer 
to $g(k(A)) \lesssim_\delta q$ is $i_s$, which also 
correctly answers $f(A) \lesssim_\eps p$, as desired.
\end{proof}

The initial segment of the $\leqm$-degrees contains some 
naturally recognizable classes which are blurred together 
by the $\alpha$ rank.  The lowest $\leqm$ degree consists 
of the constant functions; right above that is the degree of the
continuous non-constant functions.  Next above that are 
two incomparable $\leqm$-degrees: the upper semi-continuous 
functions and the lower semi-continuous functions. 

\begin{prop}\label{P9}
Let $g$ be a lower semi-continuous, discontinuous function 
(for example, $g=j_1$). 
The following are equivalent for $f \in \B_1$:
\begin{enumerate}
\item $f \leqm g$
\item $f$ is lower semi-continuous.
\item For some $e$ and some parameter $Z$, $f(A) = \ell(W_e^{A\oplus Z})$, 
where $\ell$ is the representation which maps separation names 
to real numbers.
\end{enumerate}
\end{prop}
\begin{proof}
$(1 \implies 2)$.  Given $a \in \R$, we wish to show that 
$f^{-1}((a,\infty))$ 
is open.  Let $(p_i,\eps_i)_{i<\omega}$ be an infinite 
sequence of rationals such 
that $a<p_i-\eps_i$ and $\lim p_i = a$.
Let $q_i, \delta_i$ and $k_i$ witness the defining property of
$f\leqm j_1$ for each $i$.  Now suppose 
that $f(A) > a$.  For some $i$, $f(A) > p_i+\eps_i$. 
Then the only correct answer to $f(A) \lesssim_{\eps_i} p_i$ 
is 1, so it must be that $g(k_i(A)) > q_i - \delta_i$.  
The set $C:= \{ B: g(k_i(B)) > q_i - \delta_i\}$ is open 
by the lower semi-continuity of $g$, and since 1 is a 
correct answer to $f(B) \lesssim_{\eps_i} p_i$ for every $B \in C$,
we have $C \subseteq f^{-1}((p_i-\eps_i, \infty)) \subseteq f^{-1}((a,\infty))$.

$(2 \implies 3)$  Assume $Z$ is an oracle 
which lists, for each $p$, the collection of 
rational balls contained in $f^{-1}((p,\infty))$.  
To define $W_e^{A\oplus Z}(\langle p, \eps\rangle)$,
wait to see if $A$ enters $f^{-1}((p-\eps,\infty))$.  If it does, 
enumerate the bit.  The result is a separation name of $f(A)$ 
which has the additional property that it always 
answers 1 when 1 is a permissible answer.

$(3 \implies 2)$ If $f(A) = \ell(W_e^{A\oplus Z})$, 
then $f(A) > a$ if and only if for some $p,\eps$,
$a< p-\eps$ and $\langle p,\eps\rangle \in W_e^{A\oplus Z}$,
which is an open condition.

$(2 \implies 1)$ This follows from Theorem \ref{T2}
because $g$ has rank 2 and is left-sided, and $f$ 
is either discontinuous and shares these properties, or 
$f$ is continuous, in which case $f\leqm g$ by 
Proposition \ref{prop:4}.
\end{proof}

The authors observed to Kihara that if the lattice structure 
of the Baire 1 $\leqm$-degrees would continue to higher 
Baire classes in the same pattern described in
Theorem \ref{T2}, 
the $\leqm$ reducibility could be used to extend the 
definition of the $\alpha$ rank into higher Baire classes.
After seeing these results, Kihara used a different method 
to fully describe the structure of the $\leqm$-degrees beyond 
the Baire 1 functions \cite{kihara1}, and confirmed that 
the pattern does continue.

Separately and independently of this, Elekes, Kiss 
and Vidnyansky defined a generalization of the $\alpha$,
$\beta$ and $\gamma$ ranks into the higher Baire classes
\cite{ElekesKissVidnyansky2017}.  Interestingly, they were 
able to apply their extension of the $\beta$ rank 
to solve a problem in cardinal 
characteristics, but an extension of the $\alpha$ rank 
was not suitable for that problem.  It does not seem easy to 
modify our work to get a generalization of the $\beta$ rank.
We leave a more detailed discussion of the relation between 
the various generalizations to future work.

\section{A reducibility between $\leqm$ and $\leqtt$}

There is a reducibility notion which captures the $\alpha$ rank
precisely.  Consider a truth table reduction $f\leqtt g$ which 
looks at only one bit of $g$, but may use finitely many bits of $A$.

\begin{defn}
We say $f \leqtto g$ if for all rationals $p,\eps$, 
there is a continuous $k:2^\omega\rightarrow 2^\omega$, rationals 
$q,\delta$, a number $r$, and a truth table 
$h:2^{r+1} \rightarrow \{0,1\}$ such that for every 
$A \in 2^\omega$, if $b$ is a correct answer to 
$g(k(A))\lesssim_\delta q$, then $h(A\restrict r, b)$ 
is a correct answer to $f(A) \lesssim_\eps p$.
\end{defn}

\begin{prop} The relation $f\leqtto g$ is transitive.
\end{prop}
\begin{proof}
Suppose $f_1\leqtto f_2$ and $f_2 \leqtto f_3$.  Given
$p,\eps$, let $\delta,q,k,r$ and $h$ be as guaranteed 
by the fact that $f_1\leqtto f_2$.  Given $p'=q, \eps'=\delta$, 
let $k',q',\delta',r'$ and $h'$ be as 
guaranteed by the fact that $f_2 \leqtto f_3$. 
Let $r''>r'$ be also large enough that $r''$ bits of 
any input $A$ are enough to compute $r'$ bits of $k(A)$
(using compactness).  Define 
$$h''(\tau,b) = h(k(\tau)\restrict r',h'(\tau\restrict r'',b))$$
Then the reader can verify that $k'\circ k, q',\delta',
r''$ and $h''$ witness $f_1\leqtto f_3$.
\end{proof}

\begin{thm}\label{thm:tto}
If $f,g \in \B_1$, then $f\leqtto g$ if and only if $|f|_\alpha \leq |g|_\alpha$.
\end{thm}

\begin{proof}
Suppose that $f\leqtto g$.  Given $p,\eps$, let $k,q,\delta, r$ and $h$ 
witness $f\leqtto g$.  We claim that 
$\alpha(f,p,\eps) \leq \alpha(g,q,\delta)$.  The proof is very similar 
to the $\leqm$ case.  Let $Q^\nu = k^{-1}(P^\nu_{g,q,\delta})$, 
we claim that $Q^\nu$ is a derivation sequence for
$\{f^{-1}((-\infty,p+\eps)), f^{-1}((p-\eps,\infty))\}$. 
If $A \in Q^\nu\setminus Q^{\nu+1}$, then $k(A) \in P^\nu\setminus P^{n+1}$,
so for some $\tau \prec k(A)$, either $g(P^\nu \cap[\tau]) \subseteq (-\infty, q+\delta)$, or it is a subset of $(q-\delta, \infty)$.  
Without loss of generality, assume the former.  Let $\sigma\prec A$ 
be long enough that $k([\sigma]) \subseteq [\tau]$ and 
$|\sigma|\geq r$.  Then for all $A' \in [\sigma]\cap Q^\nu$, 
we have $0$ correctly answers $g(k(A'))\lesssim_\delta q$, 
and $h(\sigma\restrict r, 0)$ correctly answers $f(A) \lesssim_\eps p$. 
So $f(Q^\nu\cap[\sigma]) \subseteq (-\infty, p+\eps)$ or $(p-\eps,\infty)$.

In the other direction, suppose $|f|_\alpha \leq |g|_\alpha$.  
Since an $\leqm$ reduction is a $\leqtto$ reduction, 
Theorem \ref{T2} implies that 
it suffices to consider the successor case. Let
$\nu$ be such that $|f|_\alpha = |g|_\alpha = \nu + 1$.
It suffices to show that $f\leqtto g$ while assuming that 
$g$ is left-sided.  (The case where $g$ is right-sided is similar.)

Given $p,\eps$, let $q,\delta$ be maximal for $g$.  Let 
$\mathcal C = \{C_0,C_1\}$ and $\mathcal D =\{D_0,D_1\}$
be as in the proof of Theorem \ref{T2}.
Exactly as there, let 
$Z$ be an oracle which contains a well-order long enough
to code $\nu$, an optimal mind-change sequence $M$ for 
$\mathcal C$, and a scaffolding sequence $S$ for $\mathcal D$.

Let $r$ be long enough that $r$ bits of any input $A$ 
are enough to see when $A$ first leaves some $P^\mu_{f,p,\eps}$ 
for some $\mu \leq \nu$.  This uses compactness.

Equivalently, $r$ is long enough that for some finite initial 
segment $(\eta_j,\sigma_j,b_j)_{j<\ell}$ from $M$, 
$\cup_j [\sigma_j] = 2^\omega$, and each $|\sigma_j| \leq r$. 
Without loss of generality, we can assume that the $\sigma_j$ 
partition the space.

Define $k$ as follows.  At stage 0, 
on input $A$, let $j$ be the index for 
which $\sigma_j \prec A$.  Let $\mu_0 = \eta_j$ and $i_0=b_j$ 
and $\tau_0 = \lambda$.  Now if $b_j = 0$ (matching 
the natural left-sidedness of $g$), search through 
$S$ to find $B \in P^\nu_{g,q,\delta}\setminus D_1$, let $B_0 = B$,
and proceed exactly as in 
the proof of Theorem \ref{T2}.  But if $b_j = 1$, 
then unfortunately $P^\nu_{g,q,\delta}\setminus D_0$ is empty. 
So in this case also let $B_0 = B$ (the same one found above),
but this means $i_0$ is an incorrect answer to $g(B_0) \lesssim_\delta q$.
We will correct this later using $h$.  So if $b_j = 1$,
proceed almost exactly as in the proof of Theorem \ref{T2},
except instead of maintaining that $i_s$ is the only 
correct answer to $g(B_s) \lesssim_\delta q$, 
now maintain that $i_s$ is incorrect for that question.

The same arguments as in Theorem \ref{T2} now guarantee that 
when $\mu_s,i_s$ and $B_s$ stabilize, then 
$|A|_\mathcal C = \mu_\infty+1$, $k(A) = B_\infty$, 
and $i_\infty$ correctly answers 
$f(A) \lesssim_\eps p$.  If $b_j=0$, 
$i_\infty$ is the only correct answer to $g(B_\infty) \lesssim_\delta q$. 
If $b_j = 1$, then 
$1-i_\infty$ is the only correct
answer to $g(B_\infty) \lesssim_\delta q$.

Define $h(\sigma,b)$ as follows.  Let $j$ be the unique index 
such that $\sigma_j\prec \sigma$.  If $b_j = 0$, let $h(\sigma,b) = b$ 
(letting the doubly correct answer through).  If $b_j = 1$, 
let $h(\sigma,b) = 1-b$ (changing the only correct answer for 
$g(k(A)) \lesssim_\delta q$ into a correct answer for 
$f(A) \lesssim_\eps p$.)
\end{proof}

Pauly has alerted us that this notion is also quite natural 
in the Weihrauch framework.  Using the notation of Section
\ref{subsec:equivalent}, he asked us whether
$f\leqtto g$ if and only if $S_f \leq^c_W S_g$.
One direction is immediate; below we prove the other
using Theorem \ref{thm:tto}. 
At a first glance, the problem with going directly from 
a Weihrauch reduction to a $\leqtto$ reduction
is that a Weihrauch reduction, when 
restricted to inputs starting with $p,\eps$, 
might use several different
choices of $q,\delta$ for different parts of the domain.
A more subtle point is that in a Weihrauch reduction,
the reverse function $H$ does not need to be defined 
on all of $2^\omega\times\{0,1\}$, 
just on the collection of values 
that it could receive as input.  Therefore, we cannot 
use compactness to automatically
transform $H$ into a truth table of the kind used in
a $\leqtto$ reduction.

\begin{prop}
For all $f,g \in\mathcal B_1$, we have $f\leqtto g$ if 
and only if $S_f \leq^c_W S_g$.
\end{prop}
\begin{proof}
A $tt1$ reduction is also a Weihrauch reduction, so one 
direction is immediate.  Suppose that $S_f \leq^c_W S_g$. 
We claim that then $|f|_\alpha \leq |g|_\alpha$. 
Let $K$ and $H$ be the continuous functions witnessing the 
Weihrauch reduction.  Note that $H$ takes two arguments, 
the original input $A$, and one bit of output 
representing a correct answer to $S_g(K(A))$.
Given $p,\eps$, by compactness there are finitely many 
strings $(\sigma_i)_{i<\ell}$, and for each $i$ 
rationals $(q_{i},\delta_{i})$ such that 
$\cup_i [\sigma_i] = 2^\omega$, and $\sigma_i \prec A$ 
implies that $K(A)$ starts with $(q_i,\delta_i)$.
Let $K_1$ be defined so that 
$$K((p,\eps)\concat \sigma_i C) = (q_i,\delta_i)\concat K_1(\sigma_i C)$$
For each $i$, let $P^\nu_i = P^\nu_{g,q_i,\delta_i}$, the 
optimal derivation sequence for $g,q_i,\delta_i$. 
Define
$$Q^\nu_i = [\sigma_i] \cap K_1^{-1}(P^\nu_i),$$
and $Q^\nu = \cup_{i<\ell} Q^\nu_i$.
We claim that $Q^\nu$ is a derivation sequence for 
$\{f^{-1}((-\infty, p+\eps)), f^{-1}((p-\eps,\infty))\}$.
It suffices to check this on the restriction to each
$[\sigma_i]$ separately, as these are clopen sets.

Fix one $i<\ell$.  Suppose that 
$A \in  Q^\nu_i \setminus Q^{n+1}_i$. Then $\sigma_i \prec A$ 
and $K_1(A) \in P^\nu_i\setminus P^{\nu+1}_i$.  So for some 
$\tau \prec K_1(A)$, either
$g(P^\nu_i\cap[tau]) \subseteq (-\infty, q_i+\delta_i)$,
or it is a subset of $(q_i-\delta_i, \infty)$.  Without loss 
of generality, assume the former.  Then $(A,0)$ must be in 
the domain of $H$.  Let $b = H(A,0)$.  
Let $\sigma\prec A$ be long enough 
that $H(A',0)=b$ whenever $\sigma \prec A'$, 
and long enough that $K_1([\sigma]) \subseteq [\tau]$. 
It is a matter of definition chasing to verify that 
$f(Q^\nu_i\cap[\sigma]) \subseteq C_b$, 
where $C_0 = f^{-1}((-\infty,p+\eps))$ and 
$C_1 = f^{-1}((p-\eps, \infty))$.  This shows that $Q^\nu_i$ 
is a derivation sequence on $[\sigma_i]$, and 
thus $Q^\nu$ is a derivation sequence. 

It follows that $\alpha(f,p,\eps) \leq \max_{i<\ell} \alpha(g,q_i,\delta_i)$,
and therefore $|f|_\alpha \leq |g|_\alpha$.
\end{proof}

\section{Properties of $\leqtt$}

In this section we characterize the $\leqtt$ degrees inside $\B_1$ 
in terms of the Bourgain rank, proving part (1) 
of Theorem \ref{thm:T2}.  Define a coarsening of the order on the
ordinals as follows:
\begin{defn}
Let $\alpha \lesssim \beta$ if for every $\gamma < \alpha$, there is 
$\delta < \beta$ and $n \in \omega$ such that $\gamma < \delta \cdot n$.
\end{defn}

This coarsening is quite robust. Recall Cantor normal form for ordinals: every ordinal  $\alpha$ 
can be written uniquely as a sum of the form 
$\alpha = \omega^{\eta_1}\cdot k_1 + \dots + \omega^{\eta_n}\cdot k_n$,
where $\eta_1 > \dots >\eta_n$ and $k_i \in \mathbb N^+$. 
Considering the existence of Cantor normal form, one can see that 
$\alpha \lesssim \beta$ if for all $\eta$, $\beta \leq \omega^\eta$
implies $\alpha \leq \omega^\eta$.

The natural sum $\alpha \# \beta$ is defined by 
$\alpha \# \beta = \omega^{\xi_1}\cdot k_1 + \dots \omega^{\xi_r}\cdot k_r$, 
where $\xi_1 > \dots >\xi_r$ are exactly the exponents 
in the Cantor normal forms of $\alpha$ and $\beta$, and $k_i$ is 
the sum of the coefficients of $\omega^{\xi_i}$ in $\alpha$ and 
$\beta$.  One sees also that $\alpha \lesssim \beta$ if for every $\gamma<\alpha$,
there is $\delta < \beta$ and $n \in \omega$ such 
that $$\gamma < \underbrace{\delta \# \delta \# \dots \#\delta}_n.$$

We will show that the $\leqtt$ degrees inside $\B_1$ correspond 
to functions whose ranks are equivalent according to this relation. 
The next lemma describes the length of combined derivation sequences.

\begin{lemma}\label{L1} Let $X$ be a compact metric space and
let $\mathcal C, \mathcal D \subseteq \mathcal P(X)$. 
Let $P^\nu_\mathcal C$ and $P^\nu_\mathcal D$ be the optimal 
derivation sequences for $\mathcal C$ and $\mathcal D$. 
Let $Q^\nu$ be the optimal derivation sequence for 
$$\{C\cap D : C \in \mathcal C \text{ and } D \in \mathcal D\}.$$
Then for all $\nu$ and $\mu$,
$$Q^{\nu\#\mu} \subseteq 
P^\nu_\mathcal C \cup P^\mu_\mathcal D.$$
\end{lemma}

\begin{proof}
By induction on $\nu \# \mu$.  If $\nu \# \mu = 0$, the statement is 
immediate.  Suppose the statement holds for all pairs of ordinals 
with natural sum less than $\nu\#\mu$.  Let 
$A \not\in P^\nu_\mathcal C \cup P^\mu_\mathcal D$.  
Then there are ordinals $\eta<\nu$ and $\xi<\mu$, 
a neighborhood $U$ of $A$, and sets $C \in \mathcal C$ 
and $D \in \mathcal D$ such that 
$P^\eta_\mathcal C \cap U \subseteq C$ and
$P^\xi_\mathcal D \cap U \subseteq D$.

Let $\zeta = \max(\eta \# \mu, \nu\#\xi)$.  Then 
since 
$\eta < \nu$ and $\xi<\mu$, we have $\zeta < \nu\#\mu$. 
So by induction, 
$$Q^\zeta \subseteq Q^{\eta\#\mu} \cap Q^{\nu\#\xi} \subseteq (P^\eta_{\mathcal C} \cup P^\mu_\mathcal D) \cap 
(P^\nu_\mathcal C \cup P^\xi_\mathcal D)$$
Rearranging the right hand side, we have 
$$Q^\zeta\subseteq P^\nu_{\mathcal C} \cup P^\mu_\mathcal D \cup 
(P^\eta_\mathcal C \cap P^\xi_\mathcal D).$$
Because $P^\nu_\mathcal C \cap U = P^\mu_\mathcal D \cap U = \emptyset$ 
and $P^\eta_\mathcal C \cap P^\xi_\mathcal D \cap U \subseteq C\cap D$, 
we have $Q^{\zeta+1} \cap U = \emptyset$.  So $A \not\in Q^{\nu\#\mu}$,
because $Q^{\nu\#\mu} \subseteq Q^{\zeta +1}$.
\end{proof}

The following is then immediate by induction.
\begin{lemma}\label{lem:2}
Let $X$ be a compact metric space and let $C_i \subseteq \mathcal P(X)$ 
for all $i<r$.  Let $P_i^\nu$ be the optimal derivation sequences 
for $\mathcal C_i$, and let $Q^\nu$ be the optimal derivation sequence for 
$$\{ \cap_{i<r} C_i : C_i \in \mathcal C_i\}.$$
Then for all $(\nu_i)_{i<r}$,
$$Q^{\#_{i<r} \nu_i} \subseteq \cup_{i<r} P_i^{\nu_i}.$$
\end{lemma}

\begin{thm}\label{T3}
If $f,g \in \B_1\setminus \B_0$, then $f\leqtt g$ if and only if 
$|f|_\alpha \lesssim |g|_\alpha$.
\end{thm}
\begin{proof}
Suppose $f\leqtt g$.  
Given $p,\eps$, let $(k_i,q_i,\delta_i)_{i<r}$ and $h$ be 
as in the definition of $\leqtt$.  For each $i$, define 
$$\mathcal C_i = \{k_i^{-1}(g^{-1}((-\infty, q_i+\delta_i))),
k_i^{-1}(g^{-1}((q_i-\delta_i,\infty)))\}.$$
Let $$\mathcal C = \{\cap_{i<r} C_i : C_i \in \mathcal C_i\}.$$
We claim that any derivation sequence for $\mathcal C$ 
is also a derivation sequence for 
$$\mathcal D := \{f^{-1}((-\infty, p+\eps)), f^{-1}((p-\eps, \infty))\}.$$
This follows because 
for every $\cap_{i<r} C_i \in \mathcal C$, there 
is a $\sigma \in 2^r$ such that $\sigma(i)$ correctly answers
$g(k(A))\lesssim_{\delta_i} q_i$, 
for every $i<r$ and $A \in \cap_{i<r} C_i$. 
Therefore, for each $A\in\cap_{i<r} C_i$, $h(\sigma)$ is 
a correct answer to $f(A) \lesssim_\eps p$.  Therefore, for some 
$D \in \mathcal D$, we have $\cap_{i<r} C_i \subseteq D$, 
and the claim follows by Proposition \ref{prop:7}.

Define $Q_i^\nu = k_i^{-1}(P_{g,q_i,\delta_i}^\nu)$.  By Proposition 
\ref{prop:7}, $Q_i^\nu$ is a derivation sequence for $\mathcal C_i$.
Let $\nu_i = \alpha(g,q_i,\delta_i)$, so that $Q_i^{\nu_i} = \emptyset$.
Let $Q^\nu$ be the optimal derivation sequence for $\mathcal C$. 
By Lemma \ref{lem:2},
$$Q^{\#_{i<r} \nu_i} \subseteq \cup_{i<r} Q_i^{\nu_i}.$$
Therefore, as $Q^\nu$ is also a derivation sequence for $\mathcal D$, 
we have
$$\alpha(f,p,\eps) \leq \#_{i<r} \nu_i \leq \underbrace{\nu \# \dots \# \nu}_r,$$
where $\nu = \max_i \alpha(g,q_i,\delta_i)$.
Therefore, $|f|_\alpha \lesssim |g|_\alpha$.

Now suppose that $|f|_\alpha \lesssim |g|_\alpha$.  We run
a daisy-chain of the kind of argument used 
in the $\leqtto$ case.  Given $p,\eps$, let 
$q,\delta$ and $n$ be such that 
$\alpha(f,p,\eps)< \alpha(g,q,\delta)\cdot n$,  
and $\alpha(g,q,\delta) \geq 2$. 
Letting $\nu = \alpha(g,q,\delta)$, we may also 
guarantee that
$P^{\nu-1}_{g,q,\delta} \not\subseteq (q-\delta,q+\delta)$, 
by decreasing $\delta$ if necessary.

We will define $3n$ functions $k_i$, 
all of them associated to this same pair $q,\delta$.
The functions are defined computably relative to 
an oracle which contains enough information to compute 
notations up to $\nu$ (and thus up to $\nu\cdot n$), 
a mind-change sequence $M$ for $\{f^{-1}((-\infty,p+\eps)), 
f^{-1}((p-\eps,\infty))\}$, and a scaffolding sequence 
$S$ for $\{g^{-1}((-\infty,q+\delta)),g^{-1}((q-\delta, \infty))\}$.

Fix $B_0 \in P^{\nu-1}_{g,q,\delta}$ with 
$g(B_0)\not\in (q-\delta,q+\delta)$, 
and let $b_0$ be the unique correct answer to 
$g(B_0) \lesssim_\delta q$.
Since $\nu \geq 2$, $|B_0|_{g,q,\delta}\geq 2$.

Given input $A$, the first $n$ functions $\{k_i\}_{i<n}$ 
are used to figure out in which interval 
$$I_i = [\nu \cdot i + 1,\nu\cdot (i+1)]_{i<n}$$
$|A|_{f,p,\eps}$ lies.
Define $k_i(A)$ as follows.  Copy $B_0$ until such a time as you see 
$A \not\in P^{\nu\cdot (i+1)}_{f,p,\eps}$.  If this occurs, 
switch to copying a nearby input $B_1$ with 
$|B_1|_{g,q,\delta} < |W|_{g,q,\delta}$ 
and where the unique 
correct answer to $g(B_1)\lesssim_{\delta} q$ is $1-b_0$.
That completes the description of the first $n$ functions $k_i$.
By observing
the answers for $g(k_i(A))\lesssim_\delta q$ for $i<n$, one
can determine the uniqe $i<n$ such that $A \in I_i$.

The next $n$ functions $\{k_{n+i}\}_{i<n}$ track the mind-changes of 
$f(A) \lesssim_\eps p$ under the assumption that $A \in I_i$.
Given input $A$, and letting $B_0$ and $b_0$ be as above, 
first copy $B_0$ into the output until such a time as 
you see $A \not\in P^{\nu\cdot (i+1)}_{f,p,\eps}$.  If this occurs, 
then we also know a rank $\mu_0<\nu$ and bit $i_0$ such that 
if $|A|_{f,p,\eps} = (\nu\cdot i) + \mu_0+1$, 
then $i_0$ correctly answers 
$f(A) \lesssim_\eps p$.  Let $\tau_0$ be whatever amount of $B_0$ 
has been copied so far.  Now proceed similarly as in 
Theorem \ref{T2}, but maintain the following at each stage:
\begin{enumerate}[(i)]
\item that $|A|_\mathcal C \leq (\nu\cdot i) + \mu_s+1 \leq (\nu\cdot i) + |B_s|_\mathcal D$, 
\item that the only correct 
answer to $g(B_s)\lesssim_\delta q$ is $i_s$ if $i_0 = b_0$, and the 
only correct answer is $1-i_s$ if $i_0 \neq b_0$.
\item if $|A|_\mathcal C = (\nu\cdot i) + \mu_s+1$, then $i_s$ correctly answers 
$f(A) \lesssim_\eps p$.
\end{enumerate}
Proceeding now just as in Theorem \ref{T2}, the above can 
be maintained unless $A$ leaves $P^{\nu\cdot i}_{f,p,\eps}$. In that 
case, the output of this computation will not be used, so 
one can continue to copy whatever $B_s$ is active at the moment 
this is discovered.  But if $\mu_s$, $i_s$ and $B_s$ stabilize 
to values $\mu_\infty$, $i_\infty$ and $B_\infty$, then if $A \in I_i$, 
we have $|A|_{f,p,\eps} = (\nu\cdot i) + \mu_\infty + 1$,
$k_{n+1}(A) = B_\infty$, $i_\infty$ is a correct answer to 
$f(A) \lesssim_\eps p$, and the only correct answer to 
$g(B_\infty)\lesssim_\delta q$ is either $i_s$ or $1-i_s$ 
depending on whether $i_0 = b_0$ or not.

The last $n$ functions $\{k_{2n+i}\}_{i<n}$ are simple indicator 
functions, with $k_{2n+1}$ 
copying $B_0$ and silently carrying out the same 
computation as $k_{n+i}$ until that computation finds an $i_0$ 
and a $b_0$.  If $k_{n+i}$ finds $i_0 \neq b_0$, switch 
to a nearby $B_1$ with $|B_1|_{g,q,\delta} < |B_0|_{g,q,\delta}$ 
and where the unique correct answer to $g(B_1)\lesssim_\delta q$
 is $1-b_0$.  Otherwise (including if $i_0$ is never defined), 
continue copying $B_0$.

Putting this all together, given $A$,
a truth table which has access to 
separating bits for each $g(k_i(A))$ can correctly answer
$f(A)\lesssim_\eps p$ as follows.  First use the separating 
bits of $g(k_i(A))$ for $i<n$ to find the unique $i$ 
such that $|A|_{f,p,\eps} \in I_i$.  Then query 
$g(h_{2n+i}(A))$ to 
learn whether $i_0 = b_0$ in the computation of $k_{n+i}(A)$.  
Finally, query 
$g(k_{n+i}(A))$ to obtain a bit $b$ which correctly 
answers $f(A)\lesssim_\eps p$ if $i_0 = b_0$.  If
$i_0 \neq b_0$, then $1-b$ will do for a correct answer.
\end{proof}

As a corollary we can give a short algorithmic proof of 
the following result of Kechris and Louveau, which 
is a consequence of their Lemma 5 and Theorem 8, 
and which allows them to conclude that their 
``small Baire classes'' $\mathcal B_1^\xi$ 
are Banach algebras.

\begin{cor}{\cite{KechrisLouveau1990}}
If $f, g \in \mathcal B_1$ are bounded, 
then $$|f+g|_\alpha, |fg|_\alpha \lesssim \max(|f|_\alpha,|g|_\alpha).$$
\end{cor}
\begin{proof} 
Without loss of generality we can assume that 
$|f|_\alpha \leq |g|_\alpha$, so $f \leqtt g$.  
Also, let $M\in \mathbb R$ be chosen so that all 
outputs of $f$ and $g$ lie in $[-M, M]$.

Then $f+g \leqtt g$ via the following algorithm.  
Given $A, p, \eps$, first ask finitely 
many questions of $f$ and $g$ to 
determine both $f(A)$ 
and $g(A)$ to within precision $\eps/2$ (by 
asking each function $2M/(\eps/2)$ questions of 
the form $f(A) \lesssim_{\eps/2} q_i$, where 
the $q_i$ are evenly spaced at intervals of $\eps/2$ 
in $[-M,M]$.)  Adding the two approximations 
gives an approximation to $(f+g)(A)$ which is 
correct to within $\eps$.  Use this approximation
to answer $f(A)\lesssim_\eps p$.

Similarly, $fg \leqtt g$ as follows.  Given $A,p,\eps$, 
first use finitely many questions to 
approximate $f(A)$ and $g(A)$ to within 
precision $\eps/(2M)$.  Multiplying the 
results gives an approximation
to $(fg)(A)$ that is correct to within $\eps$.
\end{proof}

\section{Further directions and open questions}\label{sec:future}

\subsection{A road not taken}
Recall that we used admissible representations to allow 
our results about functions on $2^\omega$ to extend 
to arbitrary compact separable metrizable spaces.
Another option for extending these 
reducibilities would be to transfer the definitions literally 
to the new spaces, without using representations.
For example, one could define $f\leqm' g$ to mean that 
for every $p,\eps$, there is a continuous function $k$ 
and rationals $q,\delta$ 
such that for all $x$, we have any correct 
answer to $g(k(x))\lesssim_\delta q$  is a correct answer to 
$f(x) \lesssim_\eps p$.

This option behaves very differently from the one we chose,
for if $X$ is very connected, then there are not enough 
continuous functions $k:X\rightarrow Y$ to get the same results. 
For example, we can define two left-sided, rank 3 functions 
in $\B_1([0,1])$ are not $\leqm'$-equivalent under this alternate 
definition.  Let $f_1 = \chi_{\{1/n : n \in \omega\}}$.  And let
$f_2 = \chi_S$ where 
$$ S = \{x_I^\ast : I \text{ is a middle third}\}$$
where $I$ is a middle third means that $I$ belongs to the sequence 
$(1/3,2/3)$, $(1/9,2/9), (7/9,8/9),... $ of intervals removed to 
create the Cantor set in $[0,1]$, and $x_I^\ast$ denotes the midpoint 
of $I$.

To see that $f_2 \not\leqm' f_1$ under this less robust definition 
of $\leqm'$, fix $p = 1/2$ and $\eps=1/3$; $q$ and $\delta$ will have to be 
similarly assigned since we are working with characteristic sets.  
Then any continuous $k$ that would work for the reduction would have 
to send the Cantor subset of $[0,1]$ to $0$.  For if any $z$ from 
the Cantor subset of $I$ satisfied $k(z) \in (1/(n+1), 1/n)$, 
then by pulling back $(1/(n+1), 1/n)$ via $k$, we'd find a whole 
neighborhood of $z$ mapped to $(1/(n+1), 1/n)$, impossible since 
every neighborhood of $z$ includes an element of $S$. So 
$h(1/3) = h(2/3) = 0$.  Now, what is $k(1/2)$.  It must be equal 
to $1/n$ for some $n$ or the reduction fails.  So $k([1/3,1/2])$ 
includes both $0$ and some $1/n$.  Since $k$ is continuous and 
$[1/3,1/2]$ is connected, its image is connected so also includes 
$1/m$ for all $m>n$.  But who are getting mapped to $1/m$? 
The purported reduction is wrong on $k^{-1}(1/m)$ for such $m$.

In fact $f_2$ is not even $\leqm'$ the characteristic function 
of the rationals, for a similar reason: if $k(1/3)$ is 
irrational and $k(1/2)$ is rational, then $k([1/3,1/2])$ 
contains many rationals.  

Since the characteristic function of 
the rationals is Baire 2, this alternate generalization 
produces a very different theory, which we did not pursue further.
  
\subsection{Computable reducibilities for discontinuous functions}
  
  The original motivation for this work was to devise 
  a notion of \emph{computable} reducibility between arbitrary (especially discontinuous)
  functions.  There is a well-established notion of computable 
  reducibility between continuous functions due to Miller
  \cite{Miller2004}, based on the 
  notion of computable function due to Grzegorczyk \cite{Grzegorczyk1955, Grzegorczyk1957} 
  and Lacombe \cite{Lacombe1955a,Lacombe1955b}.  A truly satisfying notion of computable 
  reducibility for arbitrary functions would have its 
  restriction to continuous functions agree with with 
  this established notion.  Unfortunately, the computable/lightface 
  versions of our reducibilities do not have this property.  The 
  reason for this, roughly speaking, is that the Weihrauch-based 
  reductions operate pointwise, whereas the established 
  computable reducibility on continuous functions makes 
  essential use of global information in the form of the modulus 
  of continuity.  Therefore, the following question remains of interest, 
  where of course satisfaction lies in the eye of the beholder.
  
  \begin{question}
  Is there a satisfying notion of computable reducibility for 
  arbitrary functions, whose restriction to the continuous functions
  is exactly continuous reducibility in the sense of Miller?
  \end{question}
    
  And of course, it would still be interesting to know more 
  about the structure of arbitrary functions under the computable 
  versions of these reducibilities.
    
  \begin{question}
  What can be said about the degree structure of $\mathcal F(X, \mathbb R)$ under the 
  computable versions of $\leqT, \leqtt$ and $\leqm$?
  \end{question}
    
    We will address further details and progress on these questions in a forthcoming paper.
    
\bibliographystyle{alpha}
\bibliography{functions.bib}

\end{document}